\documentclass[11pt,a4paper]{amsart}
\usepackage{amsmath}
\usepackage{amsthm}
\usepackage{amsfonts}
\usepackage{amssymb}
\usepackage{a4,a4wide}
\usepackage{leqno}
\usepackage[english]{babel}
\usepackage[utf8x]{inputenc}  %% les accents dans le fichier.tex
\usepackage{esint}
\usepackage{paralist}
\usepackage{graphics}
\usepackage{graphicx}
\usepackage{grffile}
\usepackage{pstricks}
\usepackage{epsfig}
\usepackage{color}
\usepackage{subfig}
\usepackage{dsfont}
\DeclareGraphicsExtensions{.eps,.pdf,.jpeg,.png}

\usepackage{graphics}
\usepackage{euscript}

\usepackage{color}
\usepackage[T1]{fontenc}
\usepackage{palatino}

\newtheorem{theorem}{Theorem}[section]
\newtheorem{lemma}[theorem]{Lemma}

\newtheorem{claim}[theorem]{Claim}
\newtheorem{proposition}[theorem]{Proposition}

\theoremstyle{definition}
\newtheorem{definition}[theorem]{Definition}

\theoremstyle{remark}
\newtheorem{remark}[theorem]{Remark}

\numberwithin{equation}{section}
\def\la{{\langle}}
\def\ra{{\rangle}}

\def\O{{\Omega}}

\def\f{{\mathcal{F}}}
\def\h{{\mathcal{H}}}
\def\I{{\mathcal{I}}}

\def\m{{\mathcal{M}}}

\def\M{{\mathcal{M}}}
\def\D{{\mathcal{D}}}

\def\K{{\mathcal{K}}}

\def\R{{\mathbb{R}}}

\def\N{{\mathbb{N}}}

\newcommand{\opm}[1]{\M[{#1}]}

\newcommand{\oph}[3]{\h_{_{{#1},{#2}}}[{#3}]}

\newcommand{\dem}[1]{\vskip 0.2\baselineskip \noindent {\bf{#1}}\vskip 0.2\baselineskip }
\newcommand{\fdem}{\vskip 0.2 pt \hfill $\square$ }

\newcommand{\nl}[2]{\|{#1}\|_{L^2{(#2)}}}
\newcommand{\nlto}[1]{\|{#1}\|_{2}}

\newcommand{\nlp}[3]{\|{#1}\|_{L^{#2}{(#3)}}}

\def\ds{\displaystyle}

\author{Olivier Bonnefon$^1$, J\'er\^ome Coville$^1$, Guillaume Legendre$^2$}

\address{$^1$INRA, UR 546 Biostatistique et Processus Spatiaux (BioSP)\\
228 route de l'Aérodrome\\
Domaine Saint Paul - Site Agroparc\\ 
F-84914 Avignon, France}
\email{olivier.bonnefon@avignon.inra.fr}
\email{jerome.coville@avignon.inra.fr}
\address{$^2$ CEREMADE (UMR CNRS 7534)\\
Université Paris-Dauphine\\
Place du Maréchal De Lattre De Tassigny\\
75775 Paris cedex 16, France}
\email{guillaume.legendre@ceremade.dauphine.fr}
\title{Concentration phenomenon in some non-local equation}
\date\today

\begin{document}
\maketitle
\centerline{{\it Dedicated to the Professor Stephen Cantrell, with all our admiration.}}
\begin{abstract}
We are interested in the long time behaviour of the positive solutions of the Cauchy problem involving the following integro-differential equation
\begin{align*}
&\partial_t u(t,x)=\left[a(x)-\int_{\O}k(x,y)u(t,y)\,dy\right]\,u(t,x)+\int_{\O}m(x,y)[u(t,y)-u(t,x)]\,dy\qquad(t,x)\in \R_+\times\O,%\\
\end{align*}
together with the initial condition $u(0,\cdot)=u_0$ in $\O$%, where the domain $\O$ is a , the functions $k$ and $m$ are non-negative  kernels satisfying integrability conditions and the function $a$ is continuous
. Such a problem is used in population dynamics models to capture the evolution of a clonal population structured with respect to a phenotypic trait. In this context, the function $u$ represents the density of individuals characterized by the trait, the domain of trait values $\O$ is a bounded subset of $\R^N$, the kernels $k$ and $m$ respectively account for the competition between individuals and the mutations occurring in every generation, and the function $a$ represents a growth rate. When the competition is independent of the trait%, that is, the kernel $k$ is independent of $x$, ($k(x,y)=k(y)$)
, we construct a positive stationary solution which belongs to %$d\mu$ in 
the space of Radon measures on $\O$. % $\mathbb{M}(\O)$.
Morever, when this ``stationary'' measure %$d\mu$
 is regular and bounded, we prove its uniqueness and show that, for any non negative initial datum in $L^{\infty}(\O)\cap L^1(\O)$, the solution of the Cauchy problem converges to this limit measure in $L^2(\O)$.  We also construct an example for which the measure is singular and non-unique, and investigate numerically the long time behaviour of the solution in such a situation. These numerical simulations seem to reveal some dependence of the limit measure with respect to the initial datum.     
\end{abstract}

\section{Introduction and Main results}\label{s:intro}
In this paper, we are interested in the evolution of a clonal population structured with respect to a phenotypic trait and essentially subjected to three processes: mutation, growth, and competition. As an example, one can think of a virus population structured by its virulence, as this trait can be easily quantified from experimental data. For such type of population, a common model used (see \cite{Burger1994,Burger2000,Dieckmann2003,Fournier2004,Calsina2005,Calsina2007,Perthame2007,Champagnat2008,Calsina2012,Raoul2011,Raoul2012}) is the following:
\begin{align}
&\partial_tu(t,x)=\opm{u}(t,x)+\left[a(x)-\int_{\O}k(x,y)u(t,y)\,dy\right]\,u(t,x),\quad (t,x)\in\R_+\times\O,\label{bcl-eq-intro}\\
&u(0,\cdot)=u_0 \quad \text{ in }\O, \label{bcl-eq-intro-ci}
\end{align}
the function $u\geq0$ being the density of individuals of the considered population characterized by the trait $x$, the set $\O$ is a bounded domain of $\R^N$, the function $k$ and $a$ respectively are a competition kernel and a growth rate, and $\m$ is a linear diffusion operator modelling the mutation process. In the literature, depending on the context, several kinds of mutation operator have been considered, see \cite{Burger1994,Burger2000,Fournier2004,Calsina2005,Barles2008,Champagnat2008,Desvillettes2008,Lorz2011,Raoul2011,Mirrahimi2013,Meleard2015} among others. In the present work, we focus our analysis on populations for which $\m$ is an integral operator of the form
\begin{equation}\label{bcl-def-m}
\forall v\in L^1(\O)\cap L^\infty(\O),\ \opm{v}(x):=\int_{\O}m(x,y)[v(y)-v(x)]\,dy,
\end{equation}
with $m$ a positive kernel satisfying some integrability conditions.

Lately, this type of equation have attracted  a lot of attention and much effort has been made in the analysis of solutions of \eqref{bcl-eq-intro}. In particular, let us mention \cite{Calsina2005,Calsina2007,Desvillettes2008,Calsina2012} for the construction of a global solution in $C^1(\R_+;L^1(\O)\cap L^{\infty}(\O))$ for any non negative initial data in $L^{\infty}(\O)$ and quite fairly general assumptions on $\O$, $k$, $m$ and $a$. We also point to \cite{Raoul2011,Calsina2012,Raoul2012} for an analysis of the existence of bounded continuous stationary solutions and their local stability for unidimensional domains $\O\subset \R$. However, the analysis of stationary solutions of \eqref{bcl-eq-intro} in higher dimension remains to be done, while the long time behaviour of positive solutions of problem \eqref{bcl-eq-intro}-\eqref{bcl-eq-intro-ci} is still not fully understood.

When mutations are neglected (that is, $m\equiv 0$), equation \eqref{bcl-eq-intro} is reduced to 
\begin{equation}\label{bcl-eq-sm}
\partial_tu(t,x)=\left[a(x)-\int_{\O}k(x,y)u(t,y)\,dy\right]\,u(t,x), \quad (t,x)\in\R_+\times\O,
\end{equation}
and, for a generic positive initial datum $u_0$, the solution to \eqref{bcl-eq-sm}-\eqref{bcl-eq-intro-ci} is known to converge weakly to a positive Radon measure $d\mu$ \cite{Diekmann2005,Desvillettes2008,Jabin2011}. This measure is, in some sense, a stationary solution of \eqref{bcl-eq-sm} representing an evolutionarily stable strategy for the system. For example, when the kernel $k$ is positive and does not depend on the trait (i.e., $k(x,y)=k(y)>0$), $d\mu$ is a measure whose support lies in the set $\Sigma:=\underset{x\in\bar\O}{\arg\max}\,a(x)$. In such a situation, one may check that a sum of Dirac masses  $d\mu=\sum_{i\in\I}\frac{a(x_i)}{k(x_i)}\,\delta_{x_i}$, with $x_i\in \Sigma$ for all $i\in\I$, is a stationary solution. When this measure $d\mu$ is unique, then the positive solution of \eqref{bcl-eq-sm}-- \eqref{bcl-eq-intro-ci} converges weakly to $d\mu$, see \cite{Jabin2011} for a detailed proof. 

Since the mutation process can be seen as a diffusion operator on the trait space, it is  expected that the long time behaviour of a positive solution to \eqref{bcl-eq-intro}--\eqref{bcl-eq-intro-ci} is simple and that such concentration phenomena does not occur. Indeed, this conjecture can be verified when the mutation operator $\m$ is a classical elliptic operator \cite{Lorz2011,Coville2012a}. When it is an integral operator as in the present situation, the existence of bounded equilibria when $\O$ is unidimensional seems to give credit to this conjecture. However, we prove that it is false in higher dimension. To this end, we exhibit a class of situations in which a positive singular measure $d\mu$, solution of \eqref{bcl-eq-intro} can be constructed, and investigate numerically the long time behaviour of positive solutions of the corresponding Cauchy problem.

\subsection{Main results}
We first state precisely the assumptions on the domain $\O$, the kernels $k$ and $m$  and the function $a$ under which the results are obtained. We suppose that the domain $\O$ is an open bounded connected set of $\R^N$ with Lipschitz boundary, that the function $a$ is such that
\begin{equation} \label{hyp1}
a\text{ is continuous over $\bar\O$ and positive},
\end{equation}
and that $m$ is a non-negative symmetric Carath\'eodory kernel function, that is, $m\ge0$, $m(x,y)=m(y,x)$ and
\begin{equation}\label{hyp2}
\forall x\in\O, m(x,\cdot)\text{ is measurable, and, for almost every $y$ in $\O$, } m(\cdot,y)\text{ is uniformly continuous.}
%&\text{ and for all compact set}\; \o\subset \O, k_{\o}(\cdot):=\sup_{x\in \o}K(x,\cdot) \in L^1(\O) \label{hyp2}
\end{equation}
Finally, we assume that the kernel $k$ is independent of the trait (i.e., $k(x,y)=k(y)$) and that it satisfies the following condition: \textit{there exist positive constants $C_0\ge c_0>0$ such that}
\begin{equation}\label{hyp3}
c_0\mathds{1}_\O\le k\le C_0\mathds{1}_\O,
\end{equation}
where $\mathds{1}_\O$ denotes the characteristic function of the set $\O$.

Let us now consider a stationary solution of \eqref{bcl-eq-intro}, that is, satisfying   
\begin{equation}\label{bcl-eq-red}
\opm{u}(x)+\left[a(x)-\int_{\O}k(y)u(y)\,dy\right]\,u(x)= 0,\quad x\in\O. 
\end{equation}

Under the above assumptions, we prove that there exists a positive Radon measure $d\mu$ solution of \eqref{bcl-eq-red} in a weak sense.

%In what follows, we denote by $\mathbb{M}(\O)$ the Banach space of Radon measures.

\begin{theorem}\label{bcl-thm1}
Assume $a,k$ and $m$ satisfy \eqref{hyp1}--\eqref{hyp3}. Then there exists a positive Radon measure $d\mu$ such that for all $\varphi$ in $C_c(\O)$, 
\begin{equation}\label{bcl-eq-redw}
\int_{\O}\left(\opm{\varphi}(x)+ a(x)\varphi(x)\right)d\mu(x)= \left(\int_{\O}\varphi(x)d\mu(x)\right)\left(\int_{\O}k(x)d\mu(x)\right).
\end{equation}
Let $\lambda_p$ be the principal eigenvalue of the operator $\m+a$ defined by
$$
\lambda_p(\m +a):= \sup\{\lambda\in\R \,|\,\exists\varphi\in C(\bar \O), \varphi>0, \; \text{s.t.}\; \opm{\varphi}+(a+\lambda)\varphi\le 0 \; \text{in} \; \O\}.
$$
Then, we have the following characterisation for the measure $d\mu$.
\begin{itemize}
\item If $\lambda_p$ is associated with an eigenfunction $\varphi_p$ which belongs to $L^1(\O)$, then  $d\mu$ is a regular (uniformly continuous) measure, that is, $d\mu=\bar u(x)dx$  with $\bar u$ in $L^1(\O)$ and is the unique strong solution of \eqref{bcl-eq-red}. Moreover, $\bar u$ is in $L^{\infty}(\O)$ when the principal eigenfunction $\varphi_p$ belongs to $L^{\infty}(\O)$.
\item Otherwise, $d\mu$ is a singular measure.   
\end{itemize}
\end{theorem}

As a consequence from the above dichotomy result, the existence of singular measure for \eqref{bcl-eq-redw} is strongly related to the non-existence of a $L^1$ eigenfunction associated with $\lambda_p$. This non-existence result has recently been established for the non-local operator $\m+a$, as shown in \cite{Coville2010,Shen2010,Coville2013c}.  

Next, we analyse the global stability of $d\mu$ and the long time behaviour of the positive solution of \eqref{bcl-eq-intro}--\eqref{bcl-eq-intro-ci}. When the measure $d\mu$ is regular, %we have a precise description of its global stability. 
we have the following result.

\begin{theorem}\label{bcl-thm2}
Assume $a$, $k$ and $m$ satisfy \eqref{hyp1}--\eqref{hyp3}. Assume there exists a positive regular Radon measure $d\mu(x)=\bar u(x)dx$ solution of \eqref{bcl-eq-red}. Assume further that $\bar u$ belongs to $L^{\infty}(\O)$. Then, for any non negative initial data $u_0\in L^{1}(\O)\cap L^{\infty}(\O)$, the positive solution $u$ to \eqref{bcl-eq-intro} -- \eqref{bcl-eq-intro-ci} satisfies
$$
\lim_{t\to \infty} \|u(t,\cdot)-\bar u\|_{L^2(\O)}=0.
$$
\end{theorem}

Note that the above global stability  implies the  uniqueness of the regular stationary positive Radon measure solution of \eqref{bcl-eq-red}. When no regular positive Radon measure exists, the convergence of a positive solution of \eqref{bcl-eq-intro} is very delicate to analyse. To shed light on the possible dynamics in such a situation, we explore numerically the behaviour of solutions of \eqref{bcl-eq-intro}--\eqref{bcl-eq-intro-ci}.

\subsection{Numerical simulations}
In order to illustrate and get some insight on the long time behaviour of solutions to \eqref{bcl-eq-intro}--\eqref{bcl-eq-intro-ci}, we numerically solve the problem for different choices of growth function $a$ and initial datum $u_0$ in two dimensions. Limiting ourselves to preliminary computations, we choose the domain as the open ball of radius $1/4$ centered at the origin, that is, $\O=B_{1/4}(0)$, and the competition and mutation kernels uniformly constant, such that $k\equiv 1$ and $m\equiv\rho$ with $\rho$ a positive constant. The system to be numerically solved thus reduces to:
\begin{align}
&\partial_t u(t,x)=\rho\left(\int_{B_{1/4}(0)}(u(t,y)-u(t,x))\,dy\right)+\left(a(x)-\int_{B_{1/4}(0)}u(t,y)\,dy\right)u(t,x),\ \text{ in }\R_+\times B_{1/4}(0),\label{bcl-eq-num}\\
&u(0,\cdot)=u_0\text{ in }\; B_{1/4}(0).
\end{align}
 
\subsubsection{A simple growth rate} 
First, we look at a situation in which the growth rate $a$ achieves its maximum at a single point, a case for which we can show the uniqueness of the stationary solution. More precisely, we have the following result.

\begin{proposition} For any positive value $\rho$, there exists a unique positive measure $d\mu$ which is a stationary solution of \eqref{bcl-eq-num}. Moreover, there exists a critical value $\rho^*$ such that the measure $d\mu$ is singular for $\rho<\rho^*$, whereas it is regular for $\rho\ge\rho^*$. In addition, for any non negative initial datum $u_0$ in $L^1(\O)\cap L^{\infty}(\O)$, the solution of \eqref{bcl-eq-num} converges weakly to $d\mu$.
\end{proposition}

\begin{figure}[!hp]
 \centering
\subfloat[$t=0$]{\includegraphics[height=2cm, width=3cm]{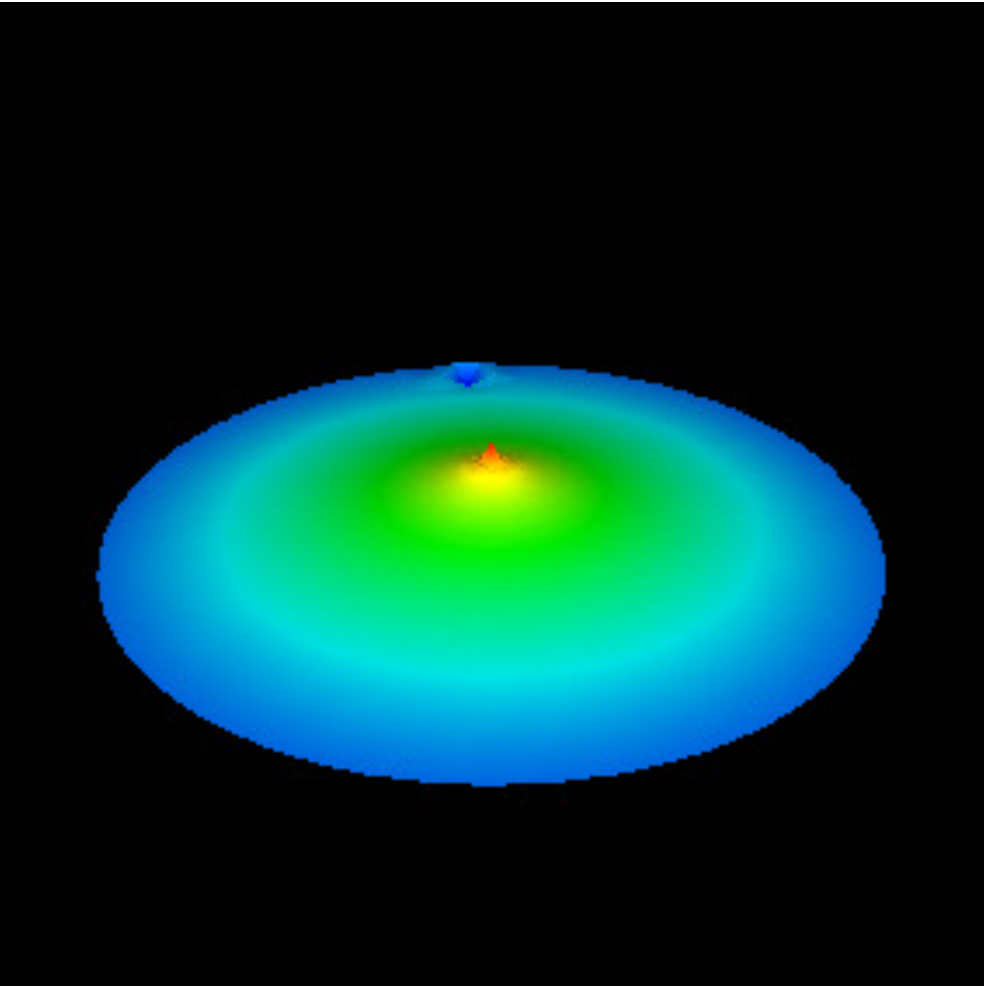}}
\subfloat[$t=10$]{\includegraphics[height=2cm,width=3cm]{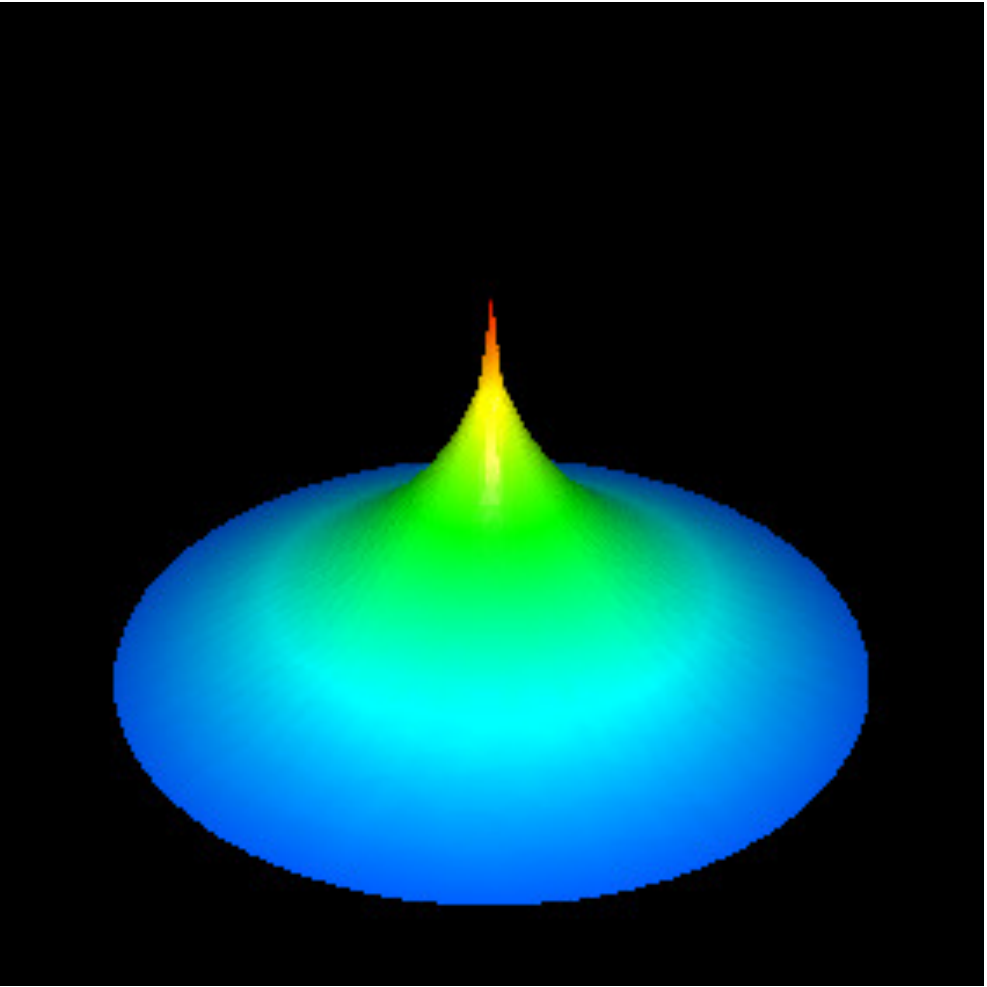}}
\subfloat[$t=20$]{\includegraphics[height=2cm,width=3cm]{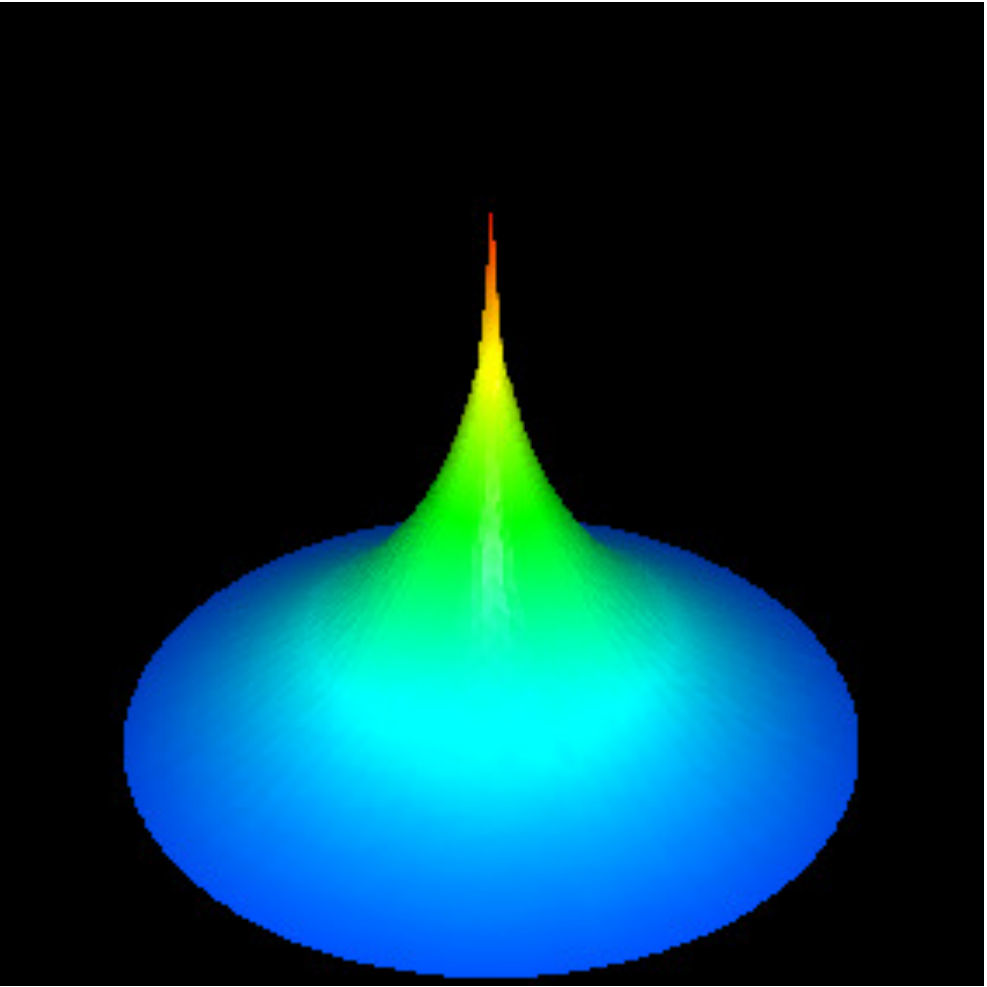}}
 \subfloat[$t=100$]{\includegraphics[height=2cm,width=3cm]{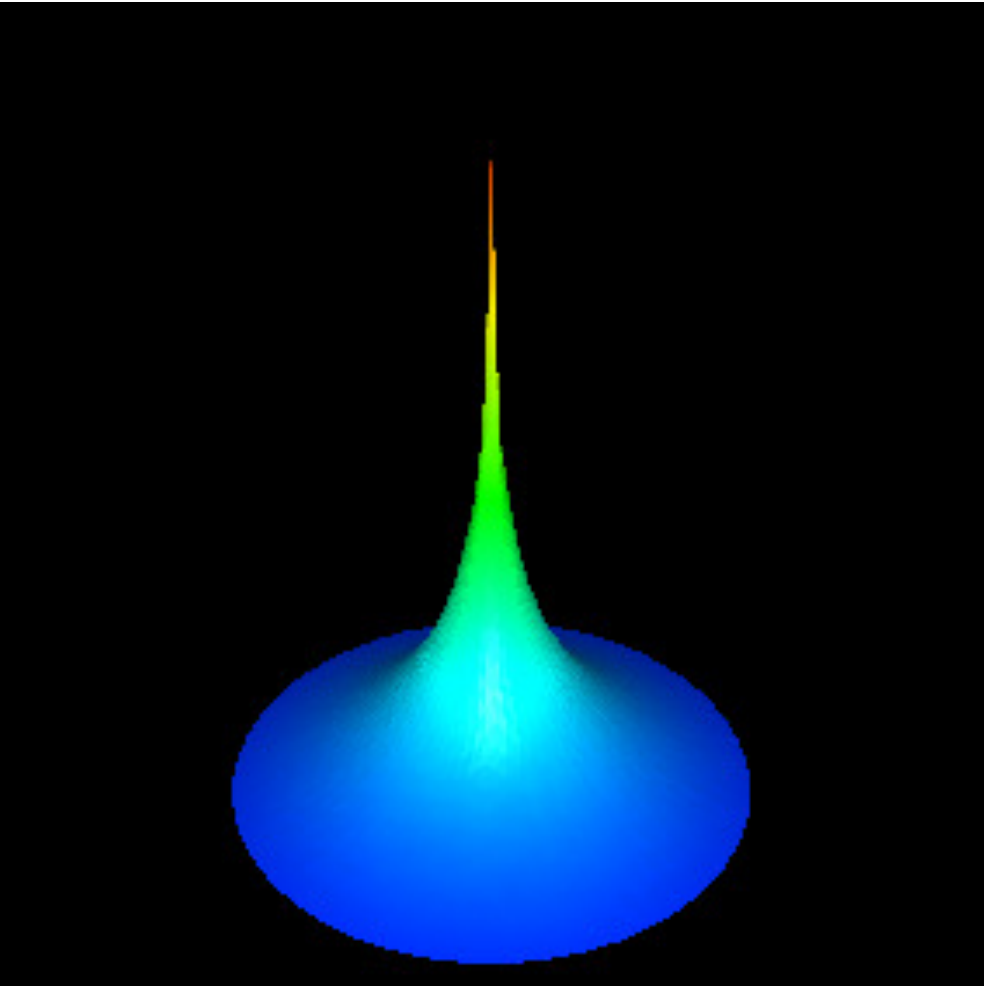}}
 %\subfloat[$t=200$]{\includegraphics[height=3cm,width=3cm]{METICE-1/DensiteT200.pdf}}
 \subfloat[$t=400$]{\includegraphics[height=2cm,width=3cm]{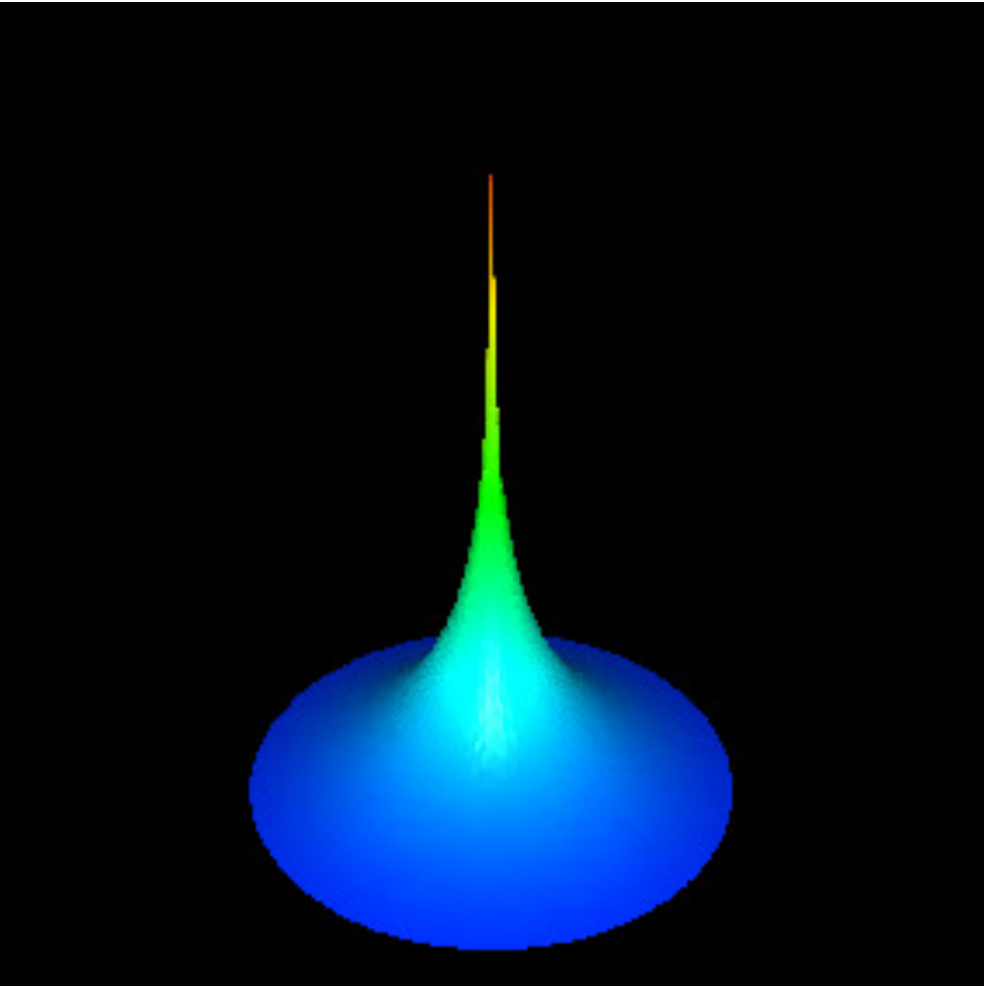}}
\qquad
\subfloat[$t=0$]{\includegraphics[width=3cm]{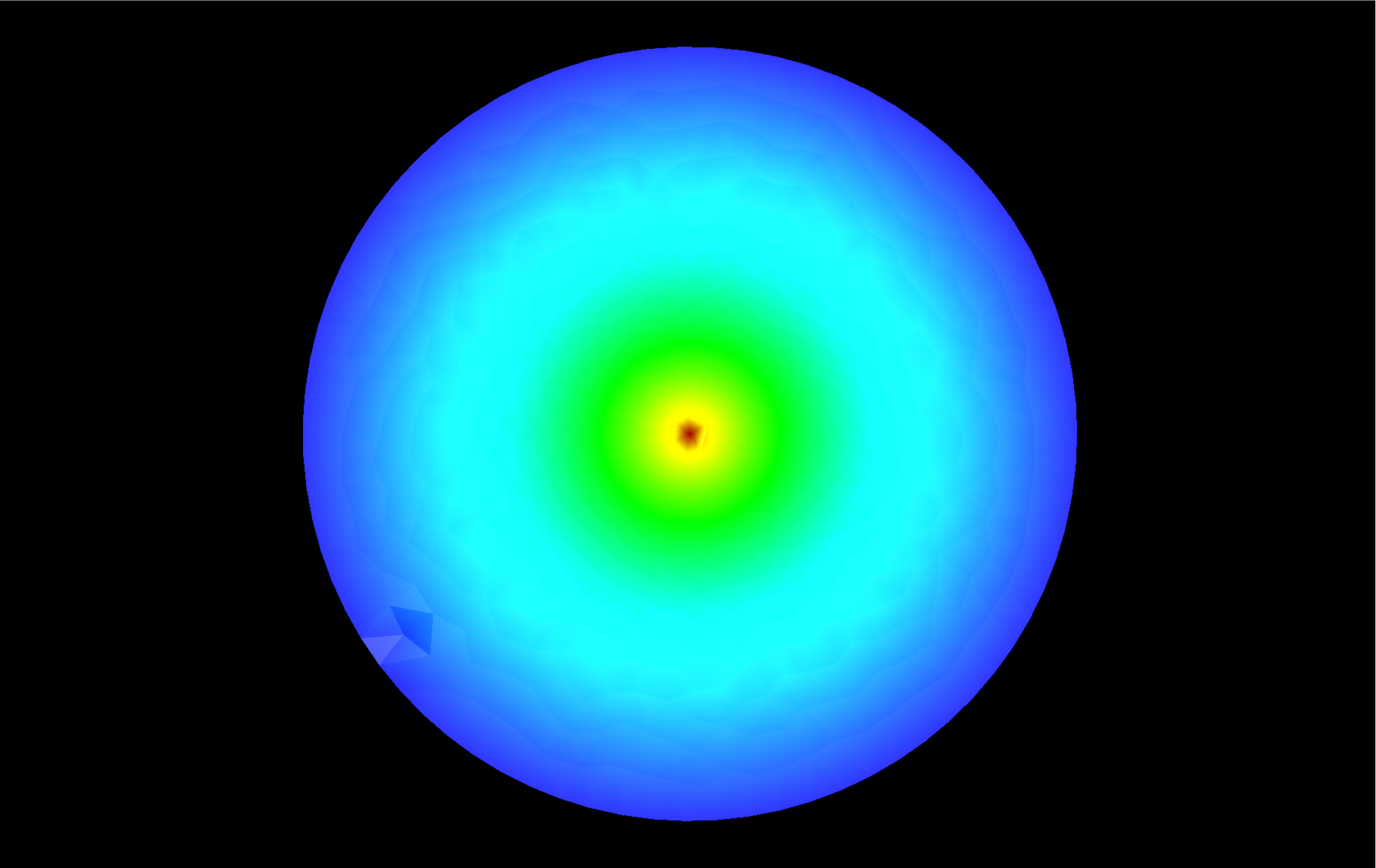}}
\subfloat[$t=10$]{\includegraphics[width=3cm]{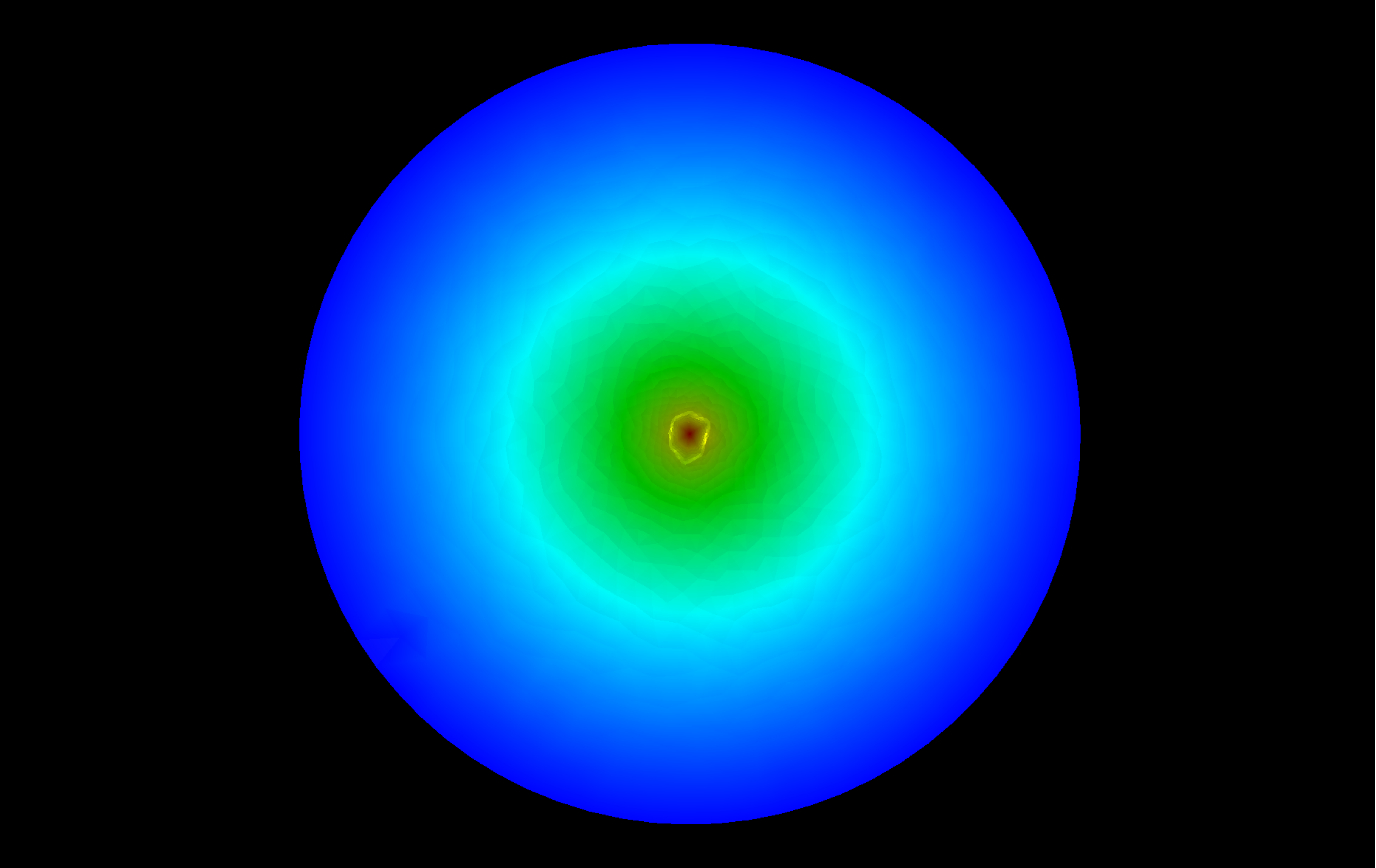}}
\subfloat[$t=20$]{\includegraphics[width=3cm]{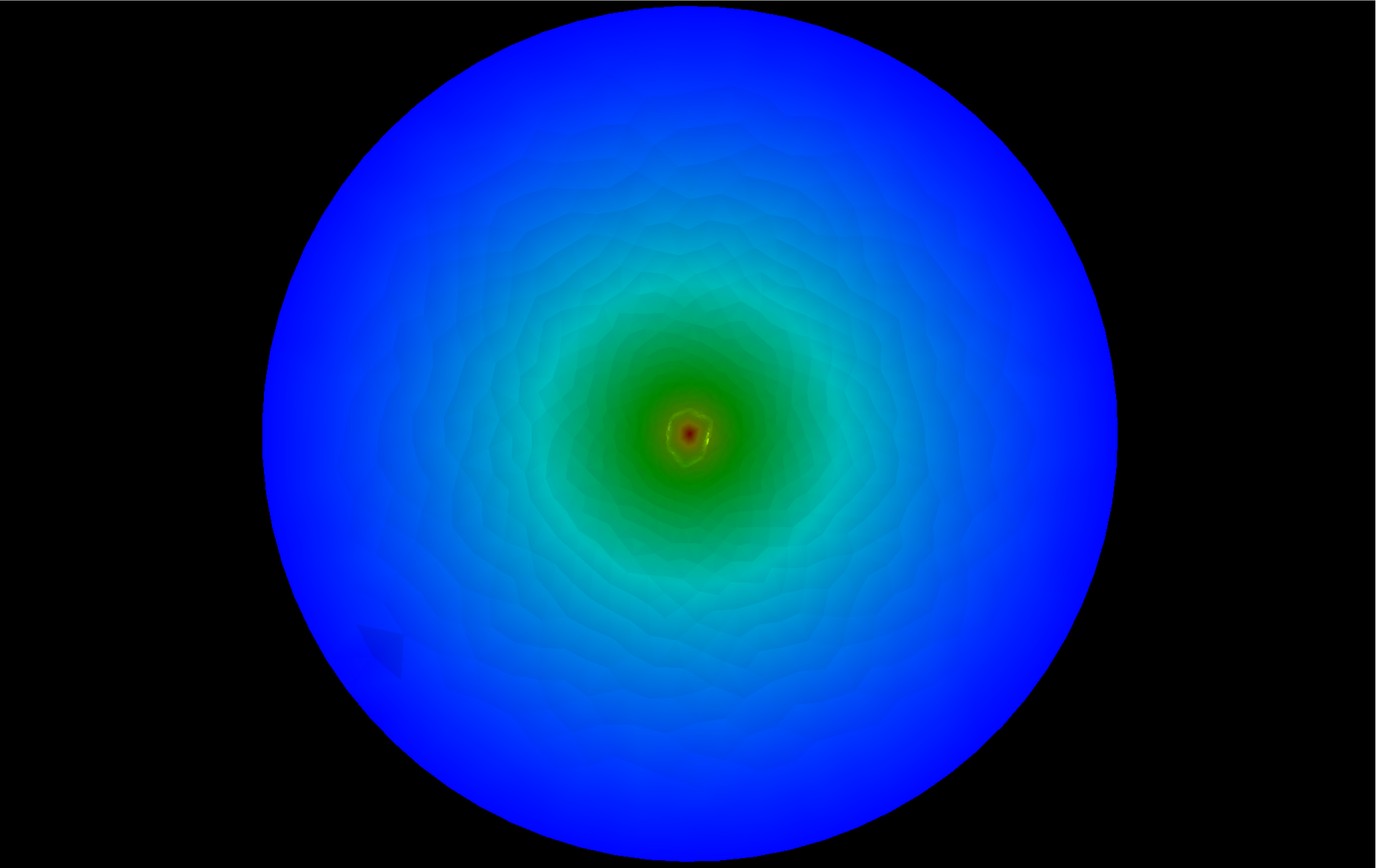}}
%\qquad
 \subfloat[$t=100$]{\includegraphics[width=3cm]{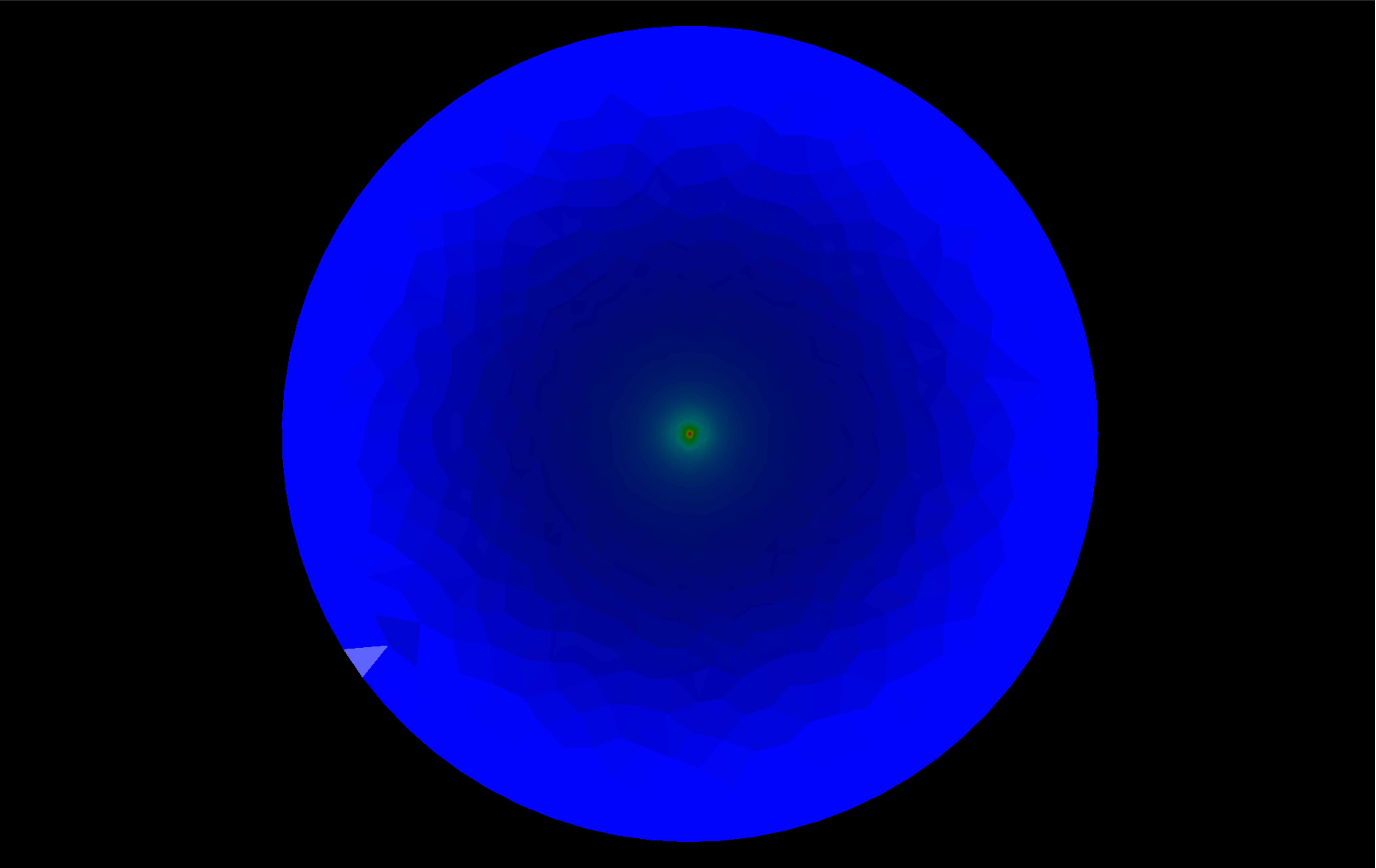}}
 %\subfloat[$t=200$]{\includegraphics[width=3cm]{images/METICE-0-1/DensiteT.200.pdf}}
 \subfloat[$t=400$]{\includegraphics[width=3cm]{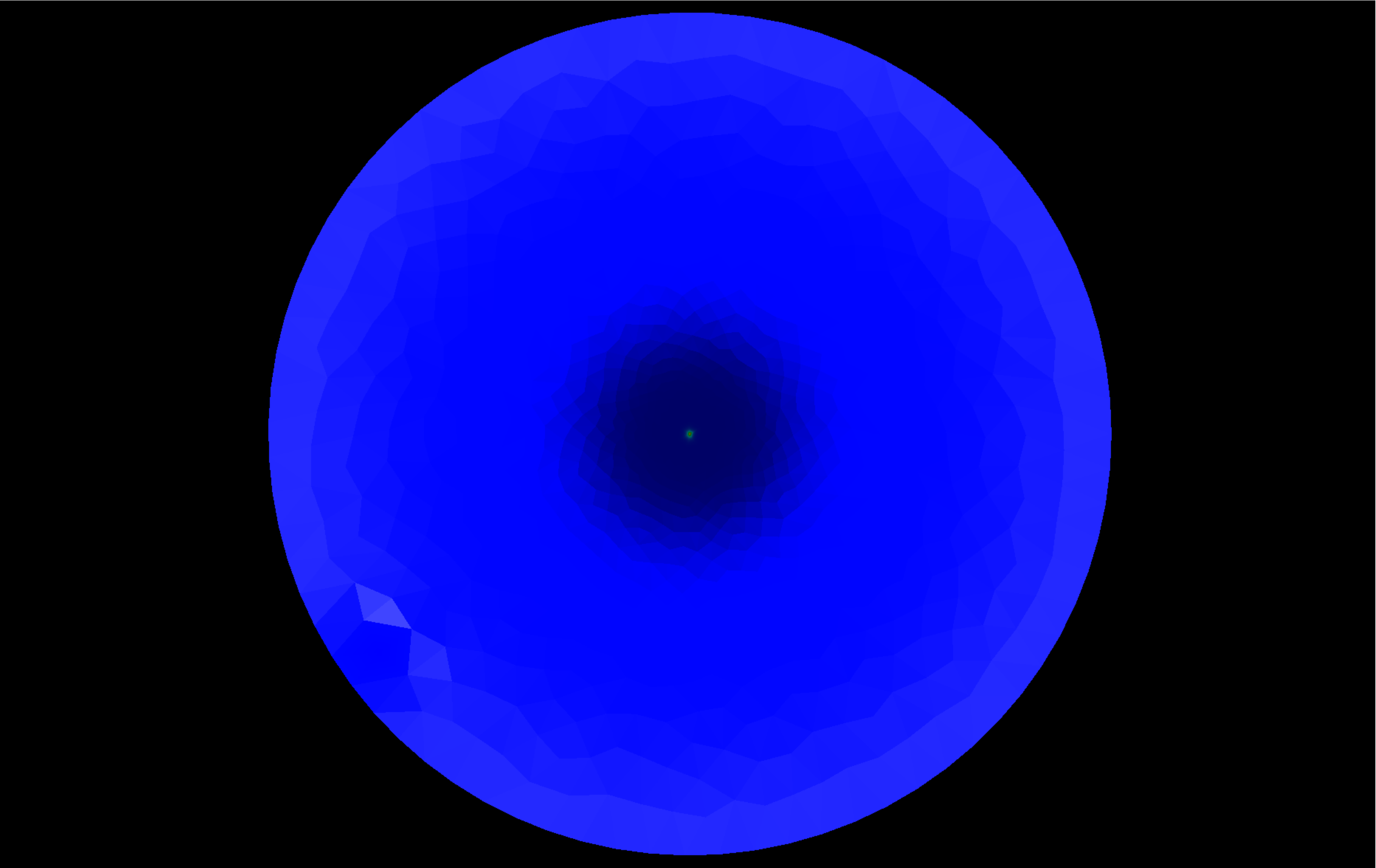}}
\caption{Numerical approximation of the solution to \eqref{bcl-eq-num} at different times for two configurations, in which the initial condition is the same and only the mutation rate differs. More precisely, we have set $\rho=1$ for the first simulation (subfigures (A) to (E)), and $\rho=0.1$ for the second one (subfigures (F) to (J)). In both situations, we observe the convergence to a stationary solution, either to a regular measure (see subfigure (E)) or to a singular measure with one Dirac mass at the origin (see subfigure (J)), the latter being characteristic of a concentration phenomenon.}\label{hdr-dg-fig-conv:1}
\end{figure}

This proposition is a direct consequence of Theorem \ref{bcl-thm1} and of the uniform $L^1$ estimates obtained in Section \ref{bcl-section-ae}. To illustrate its conclusions, we take $a(x)=1-\sqrt{\|x\|_2}$, where $\|\cdot\|_2$ denotes the Euclidean norm in $\R^2$ ($\forall x=(x_1,x_2)\in\R^2$, $\|x\|_2=\sqrt{x_1^2+x_2^2}$), and solve numerically the problem. The obtained results, presented in Figure \ref{hdr-dg-fig-conv:1}, provide a clear picture of the dynamics of the solution.

\subsubsection{A complex growth rate}
\begin{figure}[!ht]
\centering
\subfloat[$t=0$]{\includegraphics[width=3.1cm]{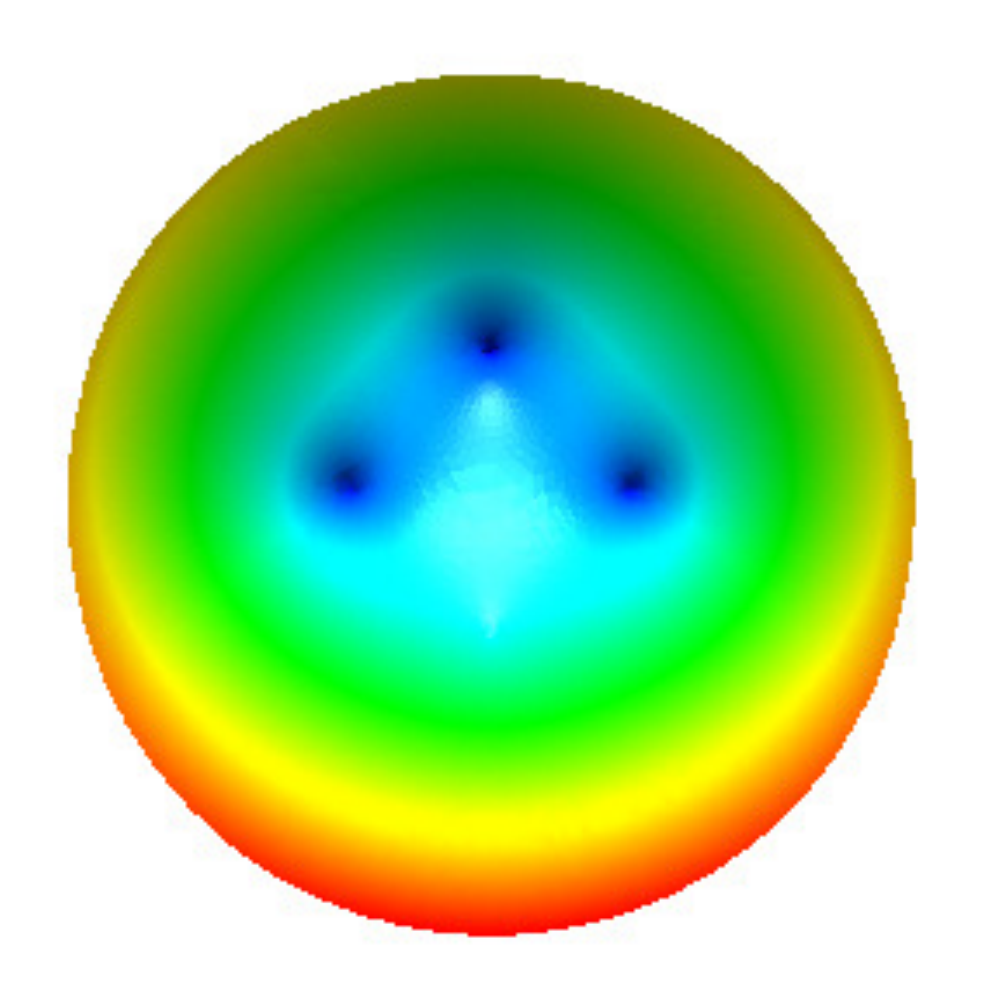}}
\subfloat[$t=50$]{\includegraphics[width=3.1cm]{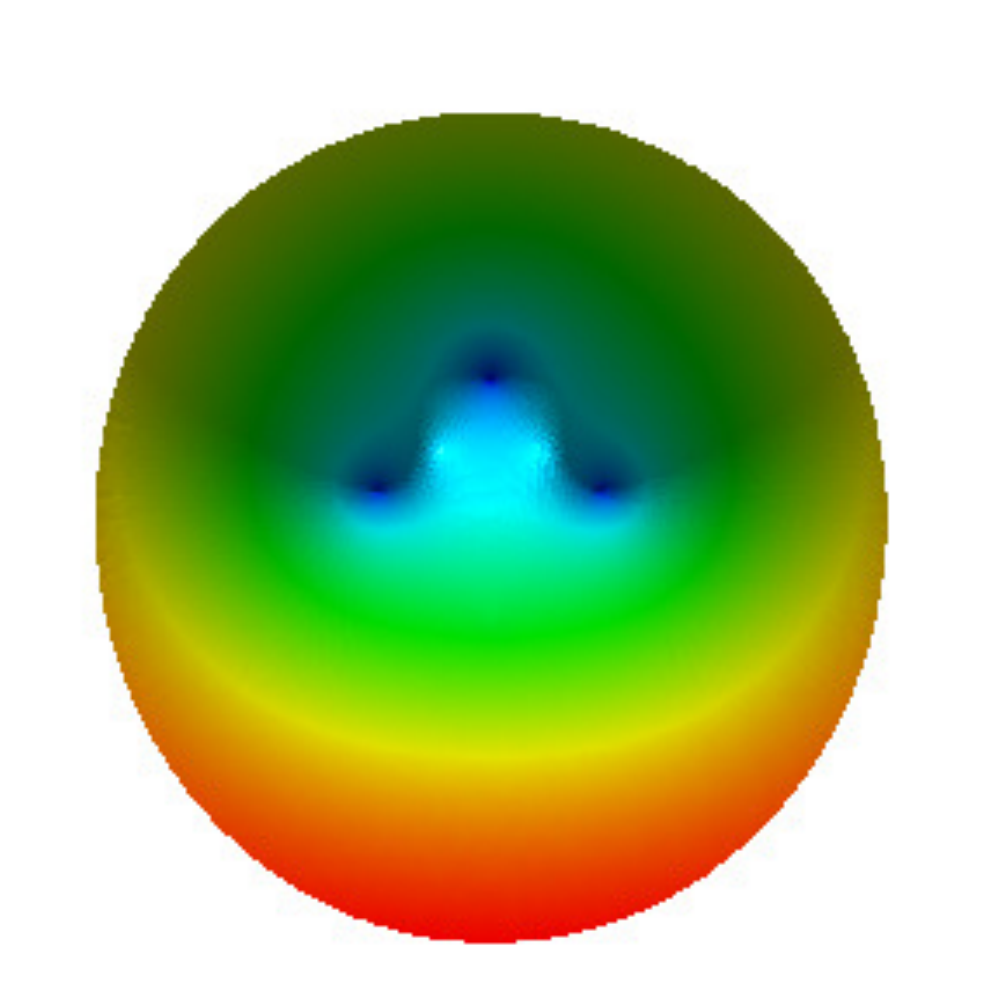}}
\subfloat[$t=100$]{\includegraphics[width=3.1cm]{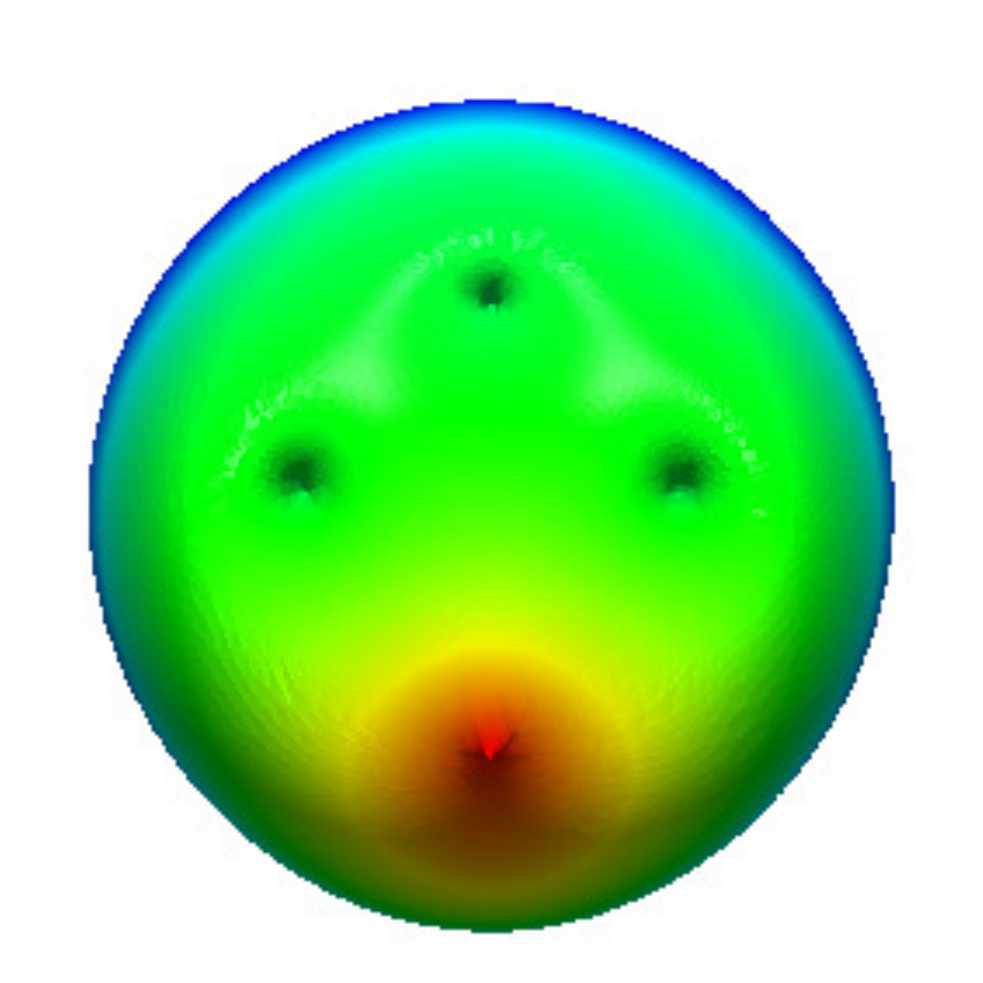}}
 %\subfloat[$t=200$]{\includegraphics[width=2.75cm]{images/rpot5-1-u0-1-img/DensiteT.200.pdf}}
 %\subfloat[$t=400$]{\includegraphics[width=3cm]{images/rpot5-1-u0-1-img/DensiteT.400.pdf}}
 %\qquad
 \subfloat[$t=800$]{\includegraphics[width=3.1cm]{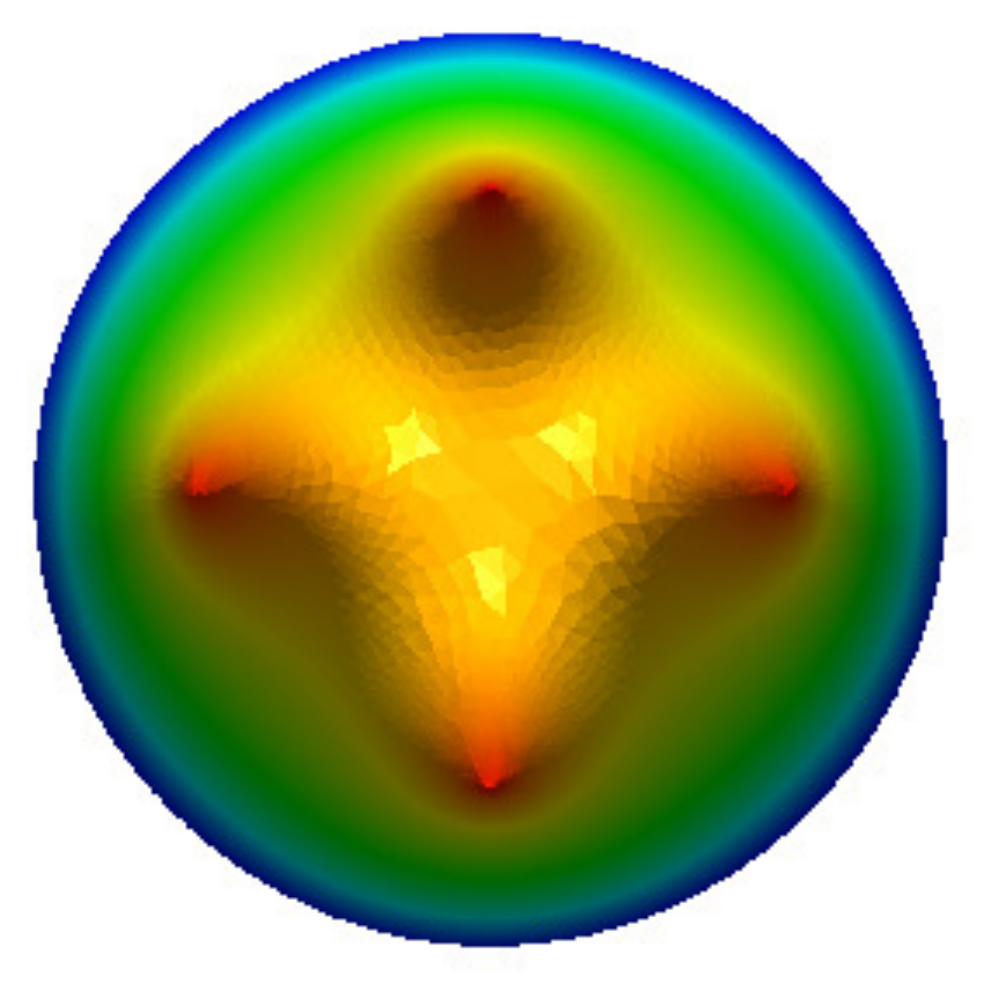}}
 \subfloat[$t=1600$]{\includegraphics[width=3.1cm]{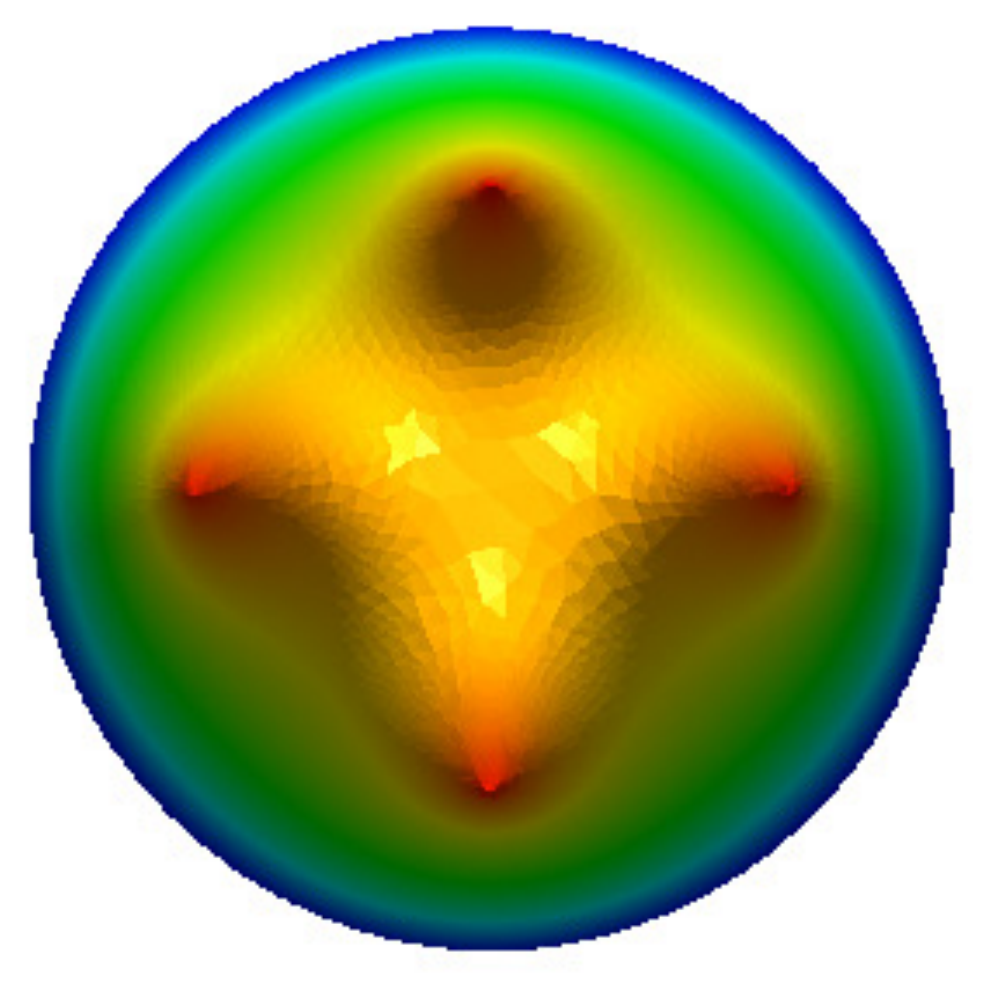}}
 \qquad
\subfloat[$t=0$]{\includegraphics[width=3.1cm]{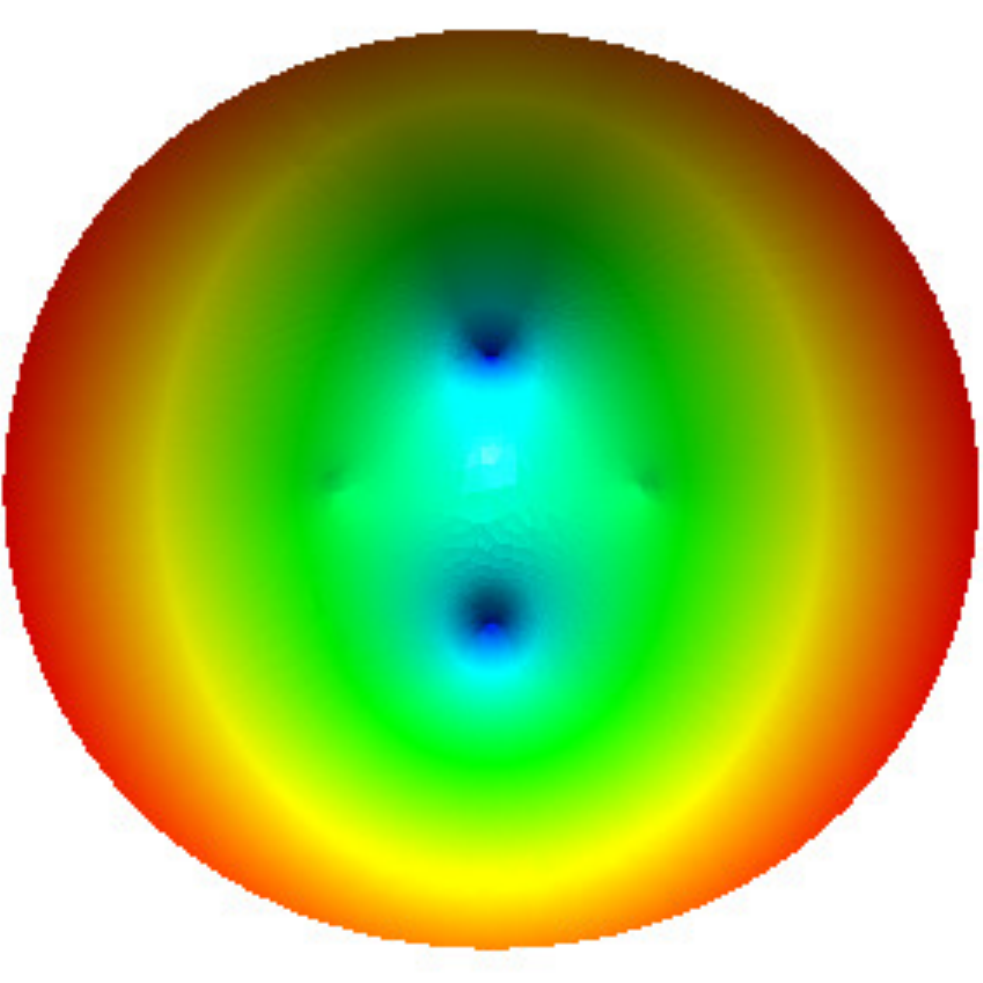}}
\subfloat[$t=50$]{\includegraphics[width=3.1cm]{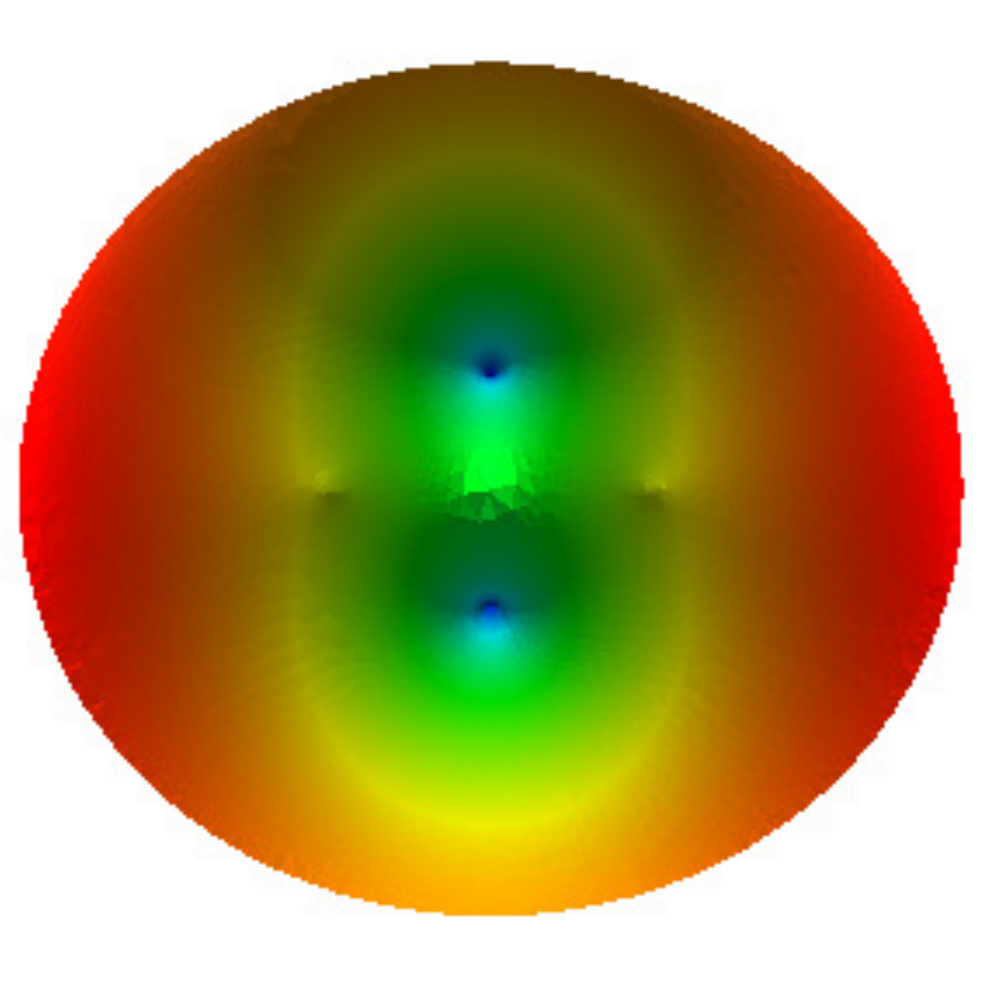}}
\subfloat[$t=100$]{\includegraphics[width=3.1cm]{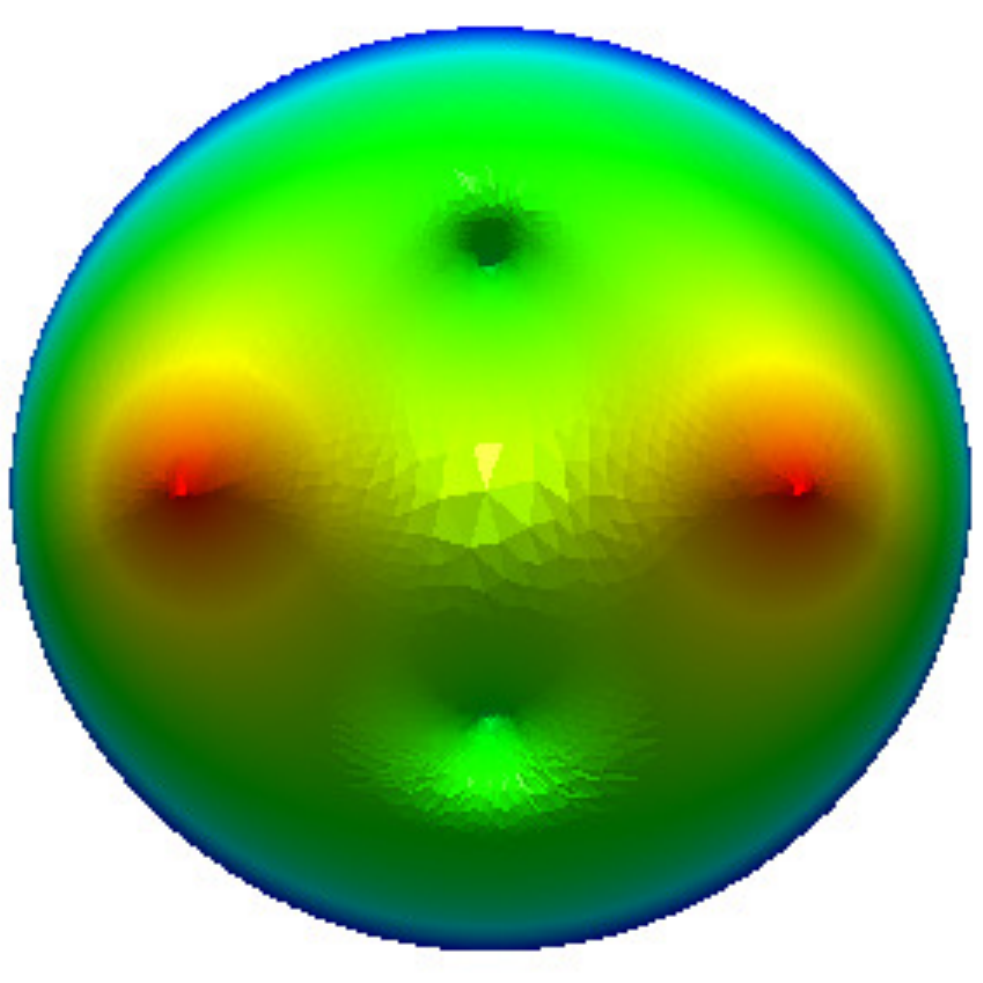}}
 %\subfloat[$t=200$]{\includegraphics[width=2.75cm]{images/rpot5-1-u0-2-img/DensiteT.200.pdf}}
 %\subfloat[$t=400$]{\includegraphics[height=2.3cm]{images/rpot5-1-u0-2-img/DensiteT.400.pdf}}
 %\qquad
 \subfloat[$t=800$]{\includegraphics[width=3.1cm]{rpot5-1-u0-2-img-DensiteT.800.pdf}}
 \subfloat[$t=1600$]{\includegraphics[width=3.1cm]{rpot5-1-u0-2-img-DensiteT.1600.pdf}}
% \qquad
%\subfloat[$t=0$]{\includegraphics[width=3.1cm]{images/rpot5-1-u0-2-img/DensiteT.1.pdf}}
%\subfloat[$t=50$]{\includegraphics[width=3.1cm]{images/rpot5-1-u0-2-img/DensiteT.50.pdf}}
%\subfloat[$t=100$]{\includegraphics[width=3.1cm]{images/rpot5-1-u0-2-img/DensiteT.100.pdf}}
% %\subfloat[$t=200$]{\includegraphics[width=2.75cm]{images/rpot5-1-u0-2-img/DensiteT.200.pdf}}
% %\subfloat[$t=400$]{\includegraphics[height=2.3cm]{images/rpot5-1-u0-2-img/DensiteT.400.pdf}}
% %\qquad
%\subfloat[$t=800$]{\includegraphics[width=3.1cm]{images/rpot5-1-u0-2-img/DensiteT.800.pdf}}
% \subfloat[$t=1600$]{\includegraphics[width=3.1cm]{images/rpot5-1-u0-2-img/DensiteT.1600.pdf}}
\caption{Numerical approximation of the solution to \eqref{bcl-eq-num} at different times for two configurations, in which the mutation rate is fixed ($\rho=1$) and the initial datum is chosen with either one (subfigures (A) to (E)) or two (subfigures (F) to (J))  spikes, that is, the function $u_0$ vanishes on the points at which the function $a$ reaches its maximum, except for one or two of them. We can see the convergence of the solution towards the same regular stationary measure (subfigures (E) and  (J)).}\label{hdr-dg-fig-conv:2}
\end{figure}

\begin{figure}[!ht]
 \centering
\subfloat[$t=0$]{\includegraphics[width=3.1cm]{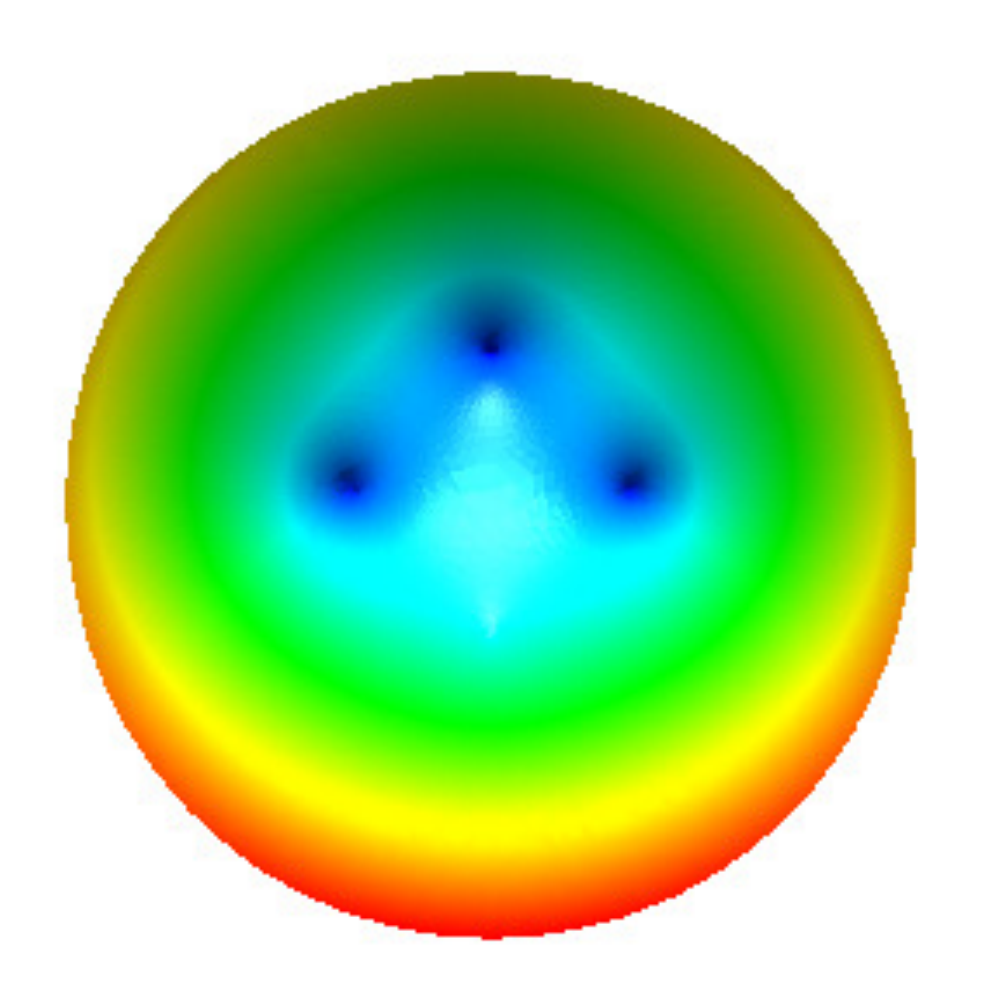}}
%\subfloat[$t=10$]{\includegraphics[width=5.5cm]{images/METICE-U0-2/DensiteT.10.pdf}}
%\subfloat[$t=20$]{\includegraphics[width=5.5cm]{images/METICE-U0-2/DensiteT.20.pdf}}
%\qquad
 %\subfloat[$t=100$]{\includegraphics[width=5cm]{images/METICE-U0-1/DensiteT.100.pdf}}
 \subfloat[$t=200$]{\includegraphics[width=3.1cm]{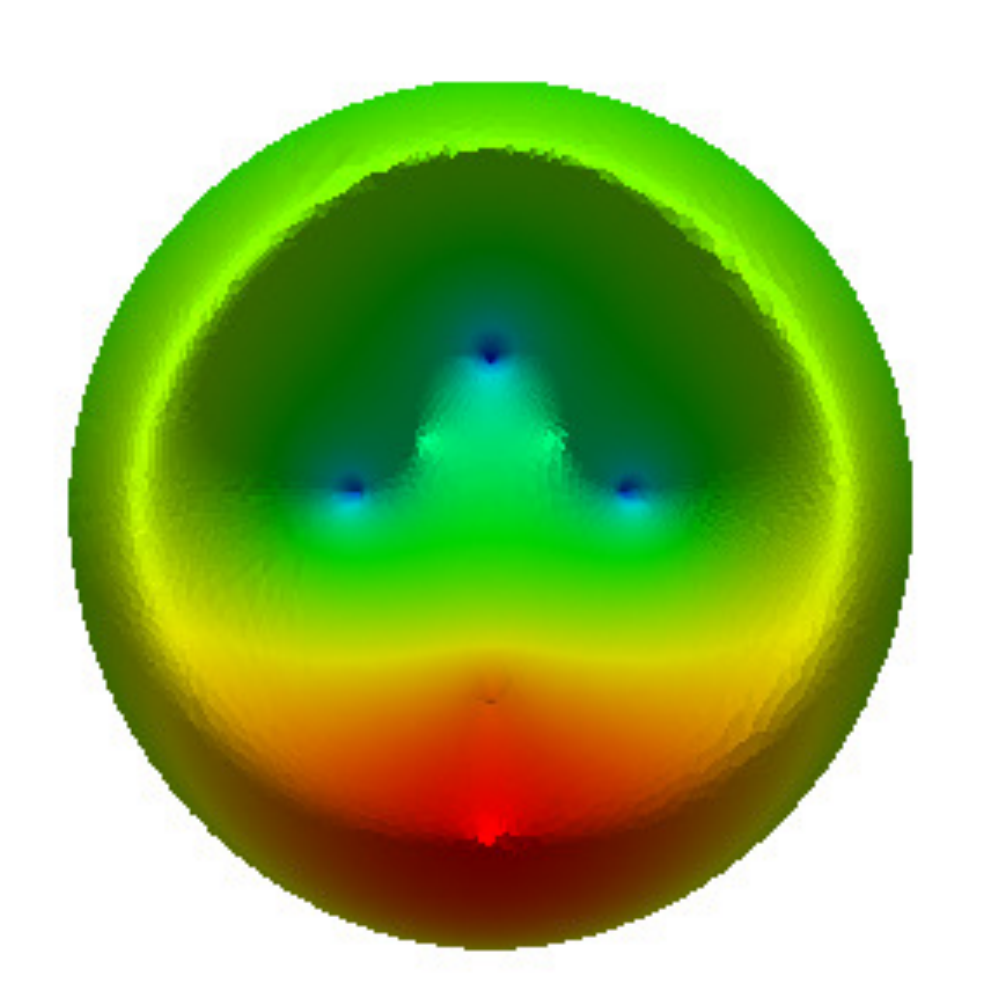}}
%\qquad 
 \subfloat[$t=400$]{\includegraphics[width=3.1cm]{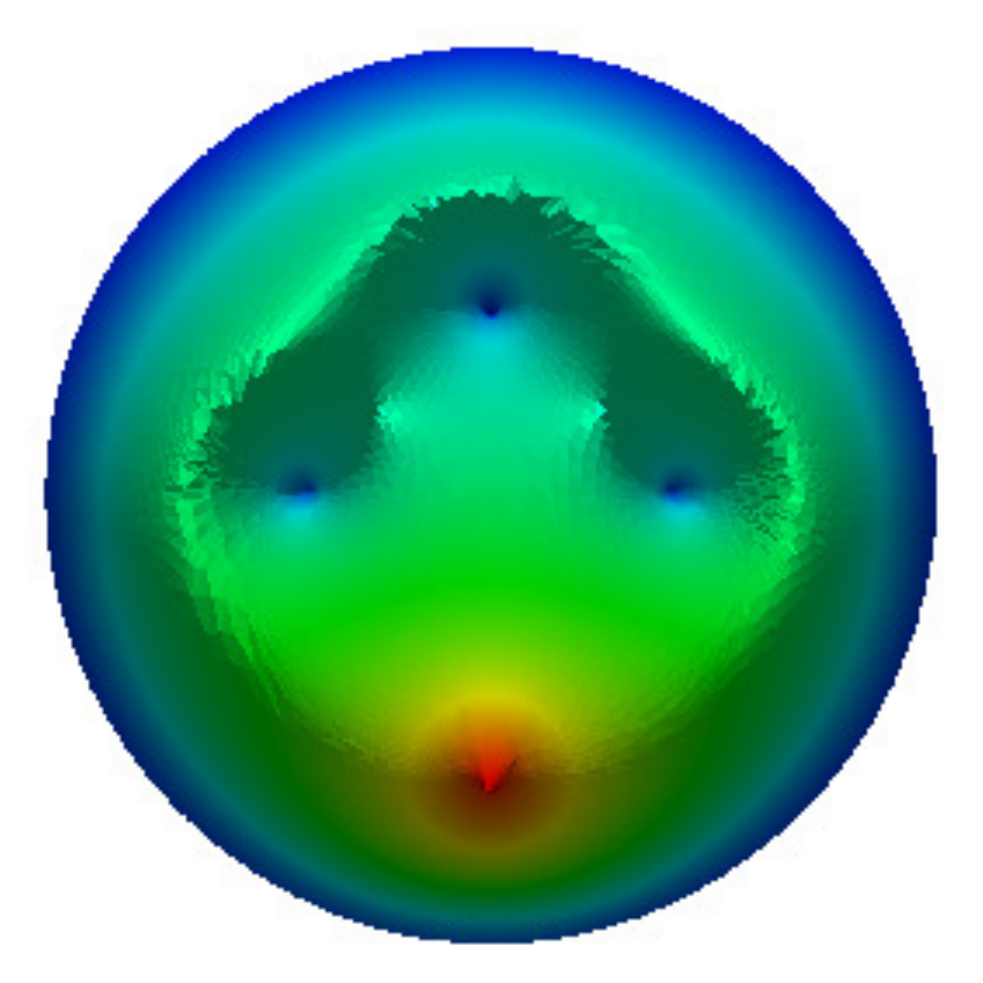}}
 %\qquad
 \subfloat[$t=800$]{\includegraphics[width=3.1cm]{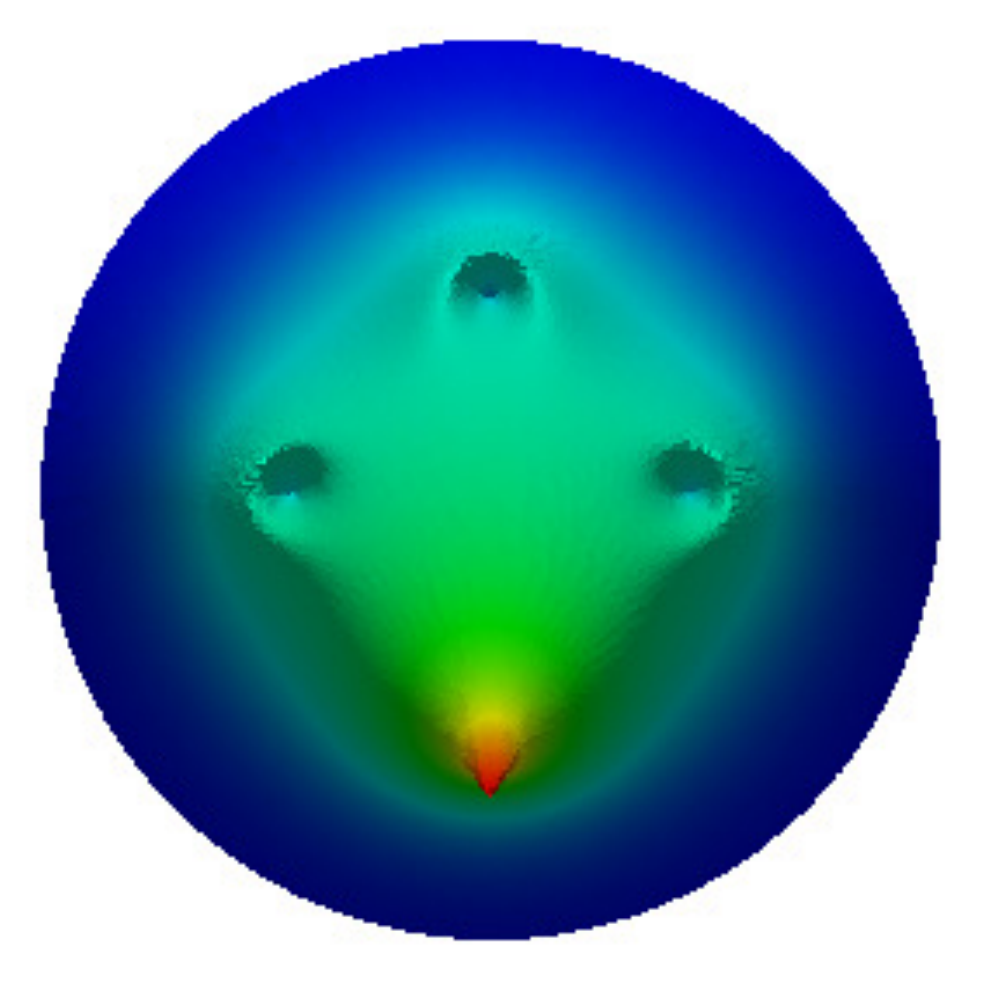}}
 \subfloat[$t=1600$]{\includegraphics[width=3.1cm]{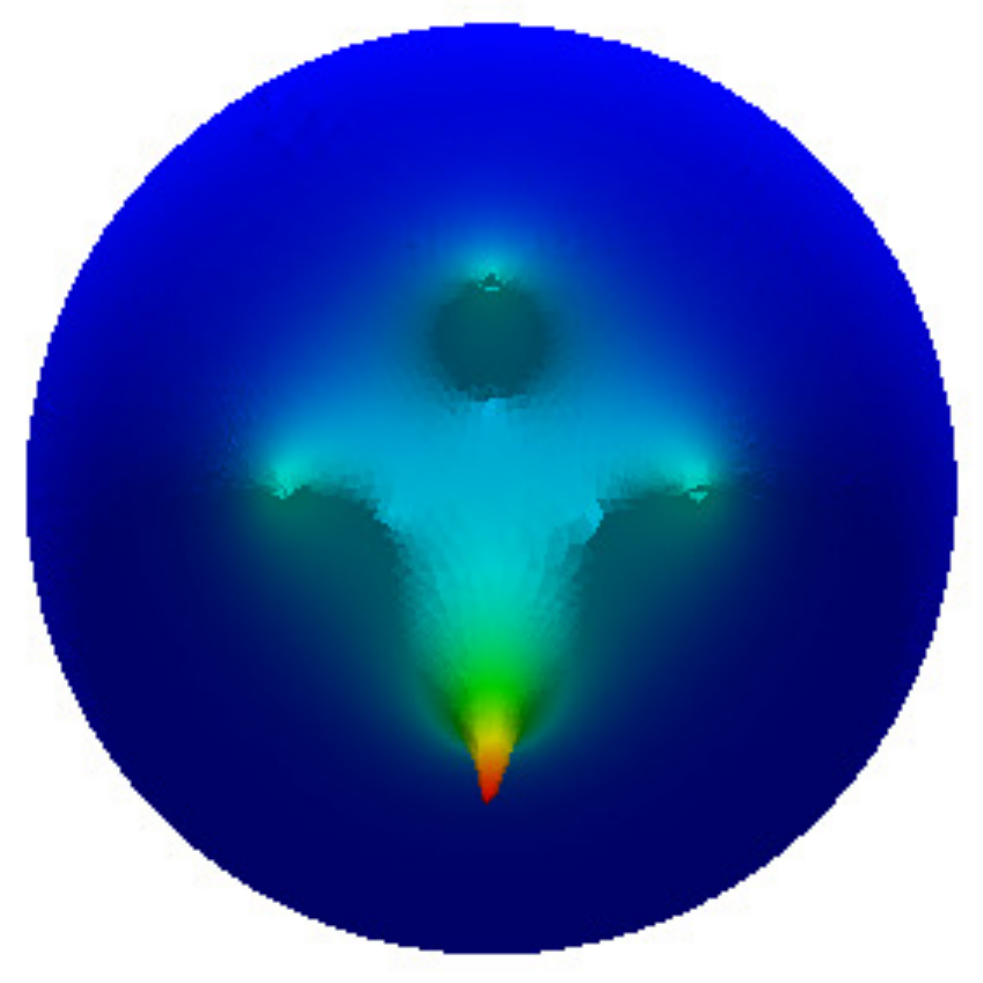}}
 %\subfloat[$t=2000$]{\includegraphics[width=5.5cm]{METICE-U0-2/DensiteT.2000.pdf}}
% \caption{Numerical approximation of the solution to \eqref{bcl-eq-num} at different times for two configurations, in which the mutation rate is fixed ($\rho=0.01$) and the initial datum is chosen with either one (subfigures (A) to (E)) or two (subfigures (F) to (J)) or three (subfigures (K) to (O)) spikes as in Figure \ref{hdr-dg-fig-conv:2}. In all the cases, we observe a concentration phenomenon: the solution converges towards a singular stationary measure, presenting one or two or three Dirac masses depending on the initial datum, which is revelatory of the non-uniqueness of the limit measure.}
\caption{Numerical approximation of the solution to \eqref{bcl-eq-num} at different times for a configuration, in which the mutation rate is fixed ($\rho=0.01$) and the initial datum is chosen with  one  spikes,  that is, the function $u_0$ vanishes on the points at which the function $a$ reaches its maximum, except for one of them. In this  cases, we observe a concentration phenomenon: the solution converges towards a singular stationary measure, presenting one Dirac mass. }%depending on the initial datum, which is revelatory of the non-uniqueness of the limit measure.} 
 \label{hdr-dg-fig-conv:3}
\end{figure}

\begin{figure}[!ht]
\subfloat[$t=0$]{\includegraphics[width=3.1cm]{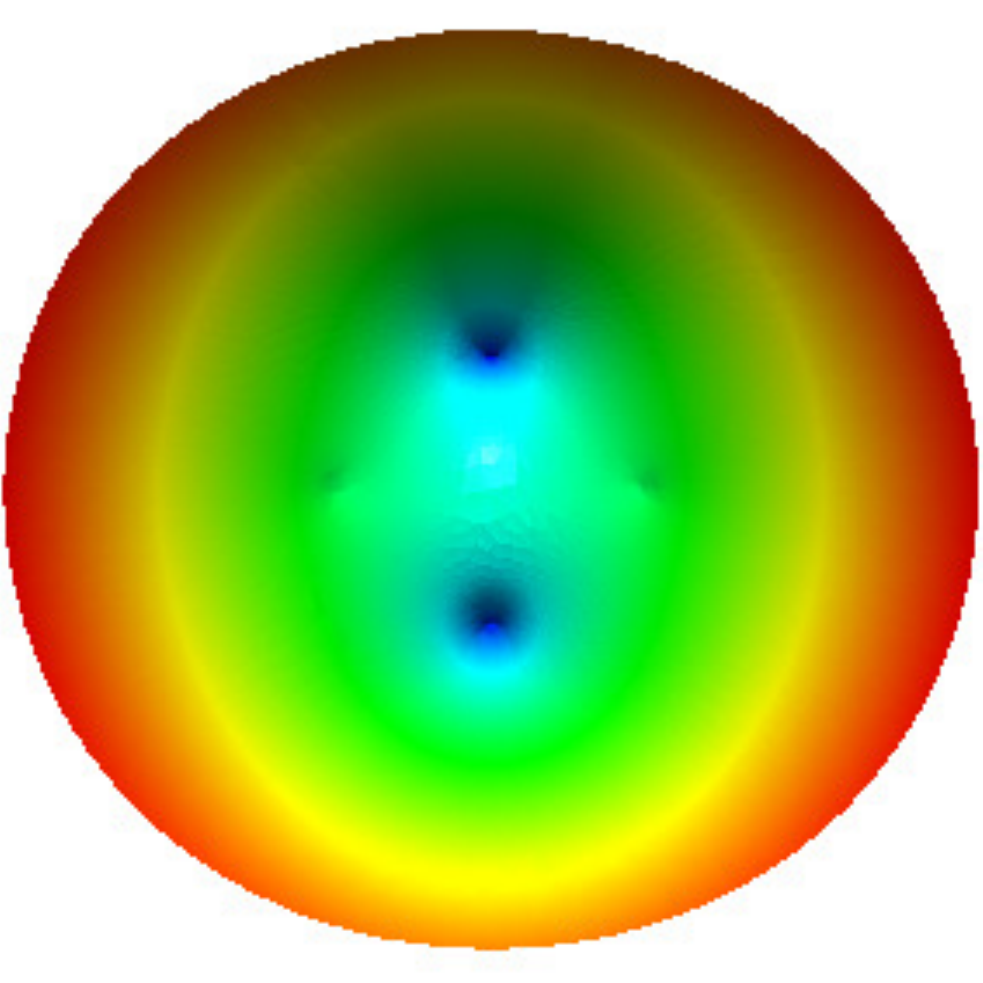}}
%\subfloat[$t=100$]{\includegraphics[width=3cm]{images/METICE-U0-2/DensiteT.100.pdf}}
 \subfloat[$t=200$]{\includegraphics[width=3.1cm]{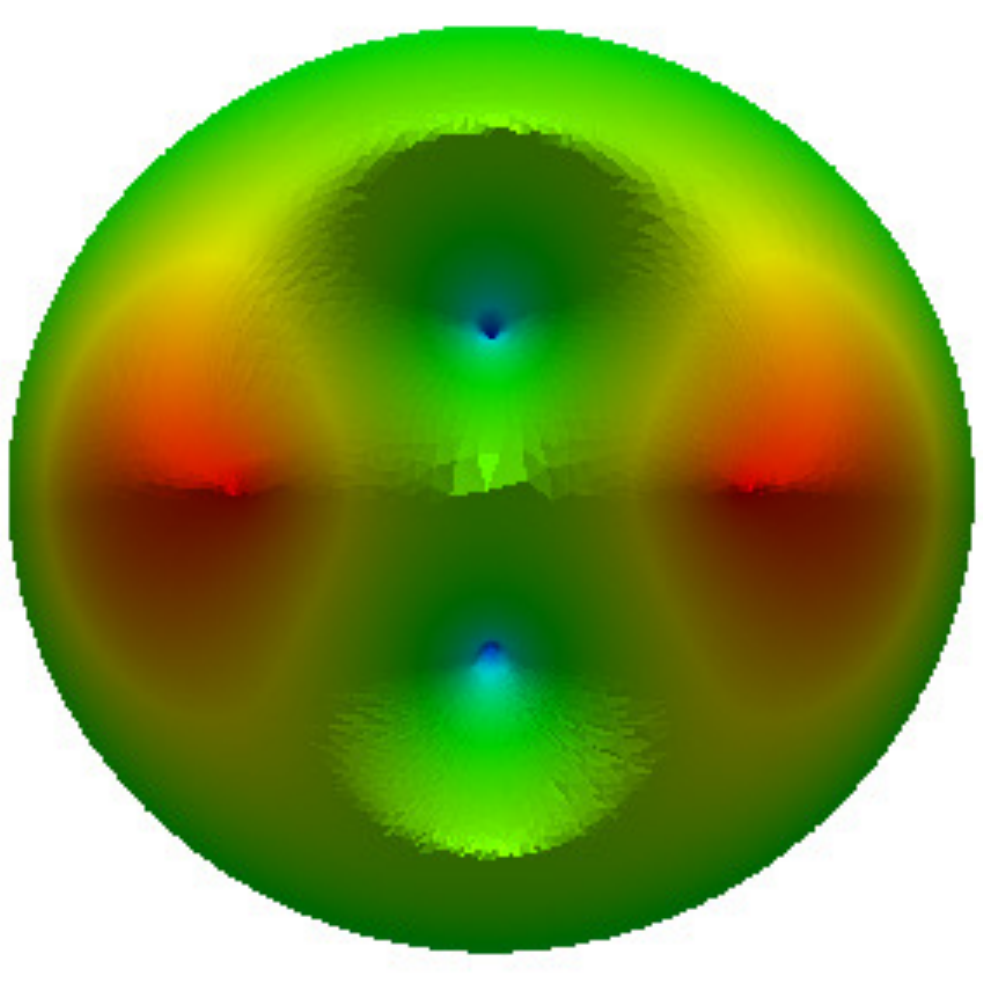}}
 \subfloat[$t=400$]{\includegraphics[width=3.1cm]{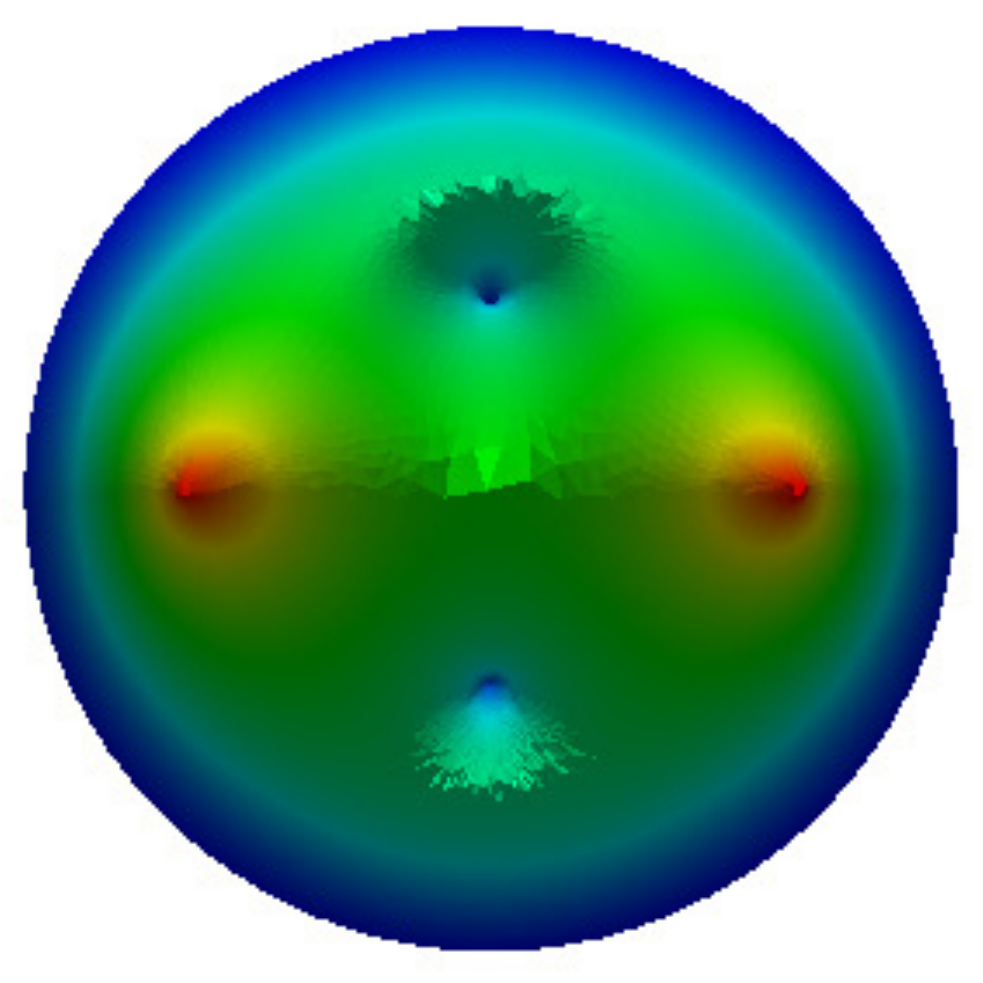}}
 %\qquad
 \subfloat[$t=800$]{\includegraphics[width=3.1cm]{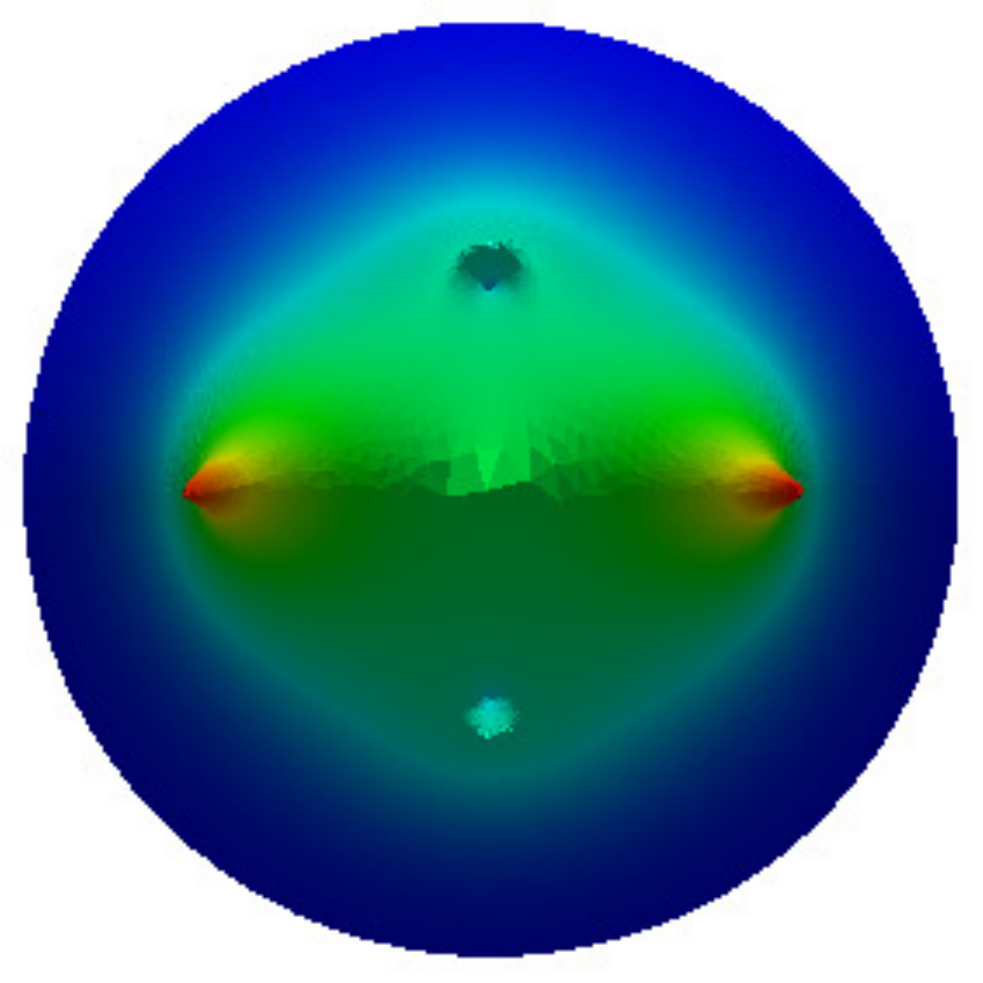}}
 \subfloat[$t=1600$]{\includegraphics[width=3.1cm]{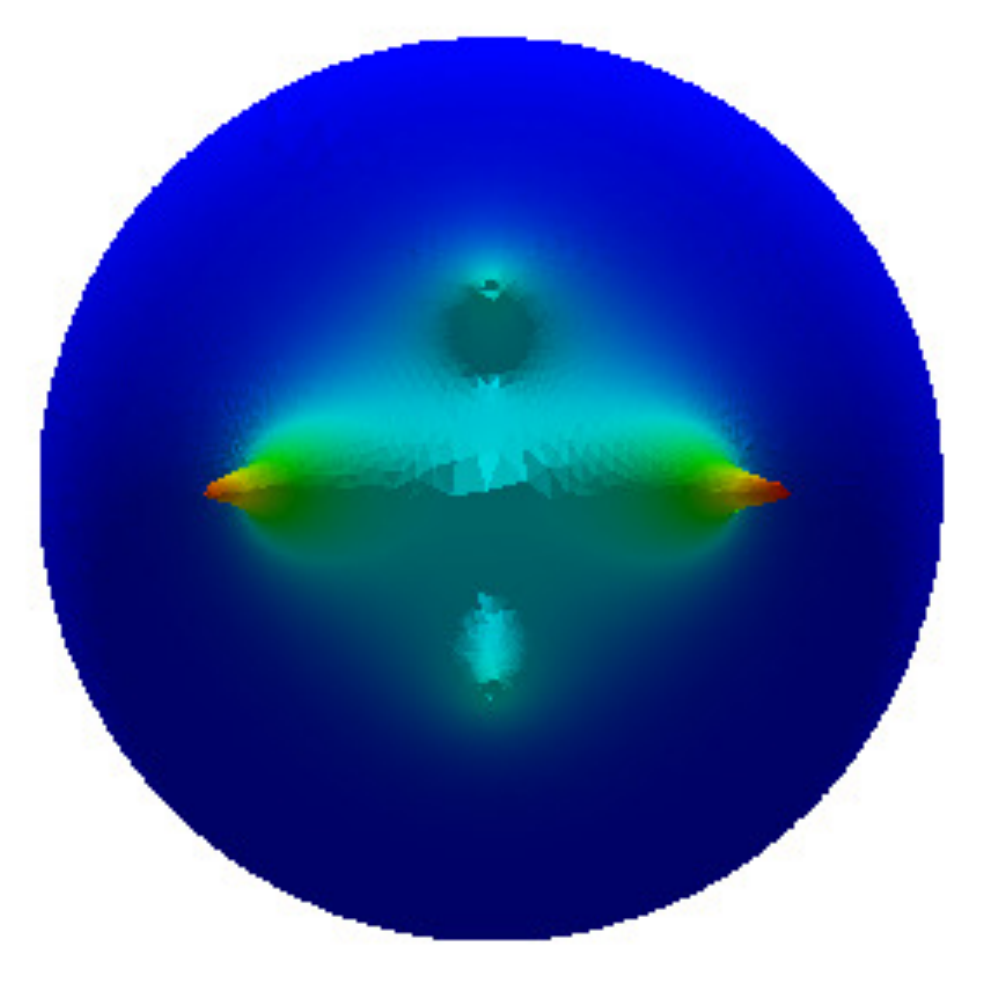}}
%\qquad
%\subfloat[$t=0$]{\includegraphics[width=3.1cm]{images/rpot5-0-01-u0-2-img/DensiteT.1.pdf}}
%%\subfloat[$t=100$]{\includegraphics[width=3cm]{images/METICE-U0-2/DensiteT.100.pdf}}
% \subfloat[$t=200$]{\includegraphics[width=3.1cm]{images/rpot5-0-01-u0-2-img/DensiteT.200.pdf}}
% \subfloat[$t=400$]{\includegraphics[width=3.1cm]{images/rpot5-0-01-u0-2-img/DensiteT.400.pdf}}
% %\qquad
% \subfloat[$t=800$]{\includegraphics[width=3.1cm]{images/rpot5-0-01-u0-2-img/DensiteT.800.pdf}}
% \subfloat[$t=1600$]{\includegraphics[width=3.1cm]{images/rpot5-0-01-u0-2-img/DensiteT.1600.pdf}} 
%\caption{Numerical approximation of the solution to \eqref{bcl-eq-num} at different times for two configurations, in which the mutation rate is fixed ($\rho=0.01$) and the initial datum is chosen with either one (subfigures (A) to (E)) or two (subfigures (F) to (J)) spikes, that is, the function $u_0$ vanishes on the points at which the function $a$ reaches its maximum, except for one or two of them. In both cases, we observe a concentration phenomenon: the solution converges towards a singular stationary measure, presenting one or two Dirac masses depending on the initial datum, which is revelatory of the non-uniqueness of the limit measure.}
\caption{Numerical approximation of the solution to \eqref{bcl-eq-num} at different times for a configuration, in which the mutation rate is fixed ($\rho=0.01$) and the initial datum is chosen with  two  spikes, that is, the function $u_0$ vanishes on the points at which the function $a$ reaches its maximum, except for two of them.  In this  cases, we observe a concentration phenomenon: the solution converges towards a singular stationary measure, presenting  two  Dirac masses.} % depending on the initial datum, which is revelatory of the non-uniqueness of the limit measure.}
\label{hdr-dg-fig-conv:4}
\end{figure}

Next, we explore a situation where the growth rate $a$ achieves its maximum at multiple points. In such a setting, we expect the stationary measure to be non-unique. In order to verify this conjecture numerically, we consider a function of the form:
$$
a(x)=1-\sqrt{\sqrt{(x_1-0.1)^2+x_2^2})}\sqrt{\sqrt{(x_1+0.1)^2+x_2^2})}\sqrt{\sqrt{x_1^2+(x_2-0.1)^2})}\sqrt{\sqrt{x_1^2+(x_2+0.1)^2})},$$  
which achieves its maximum at four distinct points. With this choice, for $\rho$ sufficiently small, we can show that there is at least four different positive Radon measures that are solution of the stationary problem \eqref{bcl-eq-redw}. The impact of the non-uniqueness of the stationary measure simulations can be seen in the simulations prensented in Figures \ref{hdr-dg-fig-conv:2} and \ref{hdr-dg-fig-conv:3}. Indeed, in a regime of mutation rate where several singular stationary measures can be constructed, we observe that the outcome of the simulation may drastically differ depending on the initial datum (see Figures \ref{hdr-dg-fig-conv:3} and \ref{hdr-dg-fig-conv:4}). In contrast, in a regime where the mutation rate is such that the stationary measure is regular, the stationary solution is a global attractor (see Figure \ref{hdr-dg-fig-conv:2}).

\subsection{Outline}
The paper is organised as follows. We start by recalling in Section \ref{bcl-section-spec} important facts about the spectral properties of the class of non-local operators considered in the problem. We then derive some uniform estimates by means of nonlinear relative entropy formulas in Section \ref{bcl-section-ae} and give a proof of Theorems \ref{bcl-thm1} and \ref{bcl-thm2} in Section \ref{bcl-section-proof}. Finally, the numerical method used for the simulations is briefly described in an appendix section. 
%the Theorem \ref{bcl-thm1}. Finally in Section \ref{bcl-section-sta} and \ref{bcl-section-asb} we prove the existence of positive steady states (Theorem \ref{bcl-thm2}) and the global stability (Theorem \ref{bcl-thm3}). A construction of a smooth positive solution to the Cauchy problem is made in the appendix.    
% \clearpage

%\end{figure}
\clearpage

\section{Spectral properties of non-local operators}\label{bcl-section-spec}
In this section, we recall some known results on the spectral problem 
\begin{equation}\label{bcv-eq-pev}
\opm{\varphi}+a\,\varphi+\lambda\,\varphi=0 \quad\text{ in }\quad \O.
\end{equation}
where $\m$ is the integral operator defined by \eqref{bcl-def-m} with a kernel satisfying the assumption \eqref{hyp2}. When the function $a$ is not constant, neither the operator $\m +a+\lambda$ nor its inverse are compact, and the Krein-Rutman theory fails in providing existence of the principal eigenvalue of $\m+a$.  However, a variational formula, introduced in \cite{Berestycki1994} to characterise the first eigenvalue of elliptic operators, %of the form $\E:=a_{ij}(x)\partial_{ij}+b_i(x)\partial_i +c(x)$,
%\begin{equation}\label{bcv-eq-lambda1-form-sup}
%\lambda_1(\E) := \sup \{\lambda \in \R \, |\, \exists \, \varphi \in W^{2,n}(\O), \varphi > 0 \;\text{ such that }\; \ope{\varphi} + \lambda\varphi\le 0 \;\text{in}\; \O\},
%\end{equation}
can be transposed to the operator $\m+a$. Namely, the following quantity
\begin{equation}\label{bcv-eq-lambda-form-sup}
\lambda_p (\m + a) := \sup \{\lambda \in \R \, |\, \exists \, \varphi \in C(\O), \varphi > 0 \;\text{ such that }\; \opm{\varphi}+a\varphi + \lambda\varphi\le 0\;\text{in}\; \O\}.
\end{equation}
is well defined and called the generalised principal eigenvalue of $\m + a$. It is known \cite{Coville2010,Shen2010,Coville2013,Berestycki2014,Shen2015} that $\lambda_p(\m+a)$ is not always an eigenvalue of $\m+a$ in a reasonable Banach space, which means there is not always a positive continuous eigenfunction associated with it. Nevertheless, as shown in \cite{Coville2013c}, there always exists an associated positive Radon measure.

\begin{theorem}[\cite{Coville2013c}]\label{bcl-thm-mu}
Let the domain $\O$ be bounded, the operator $\m$ be defined by \eqref{bcl-def-m} with a kernel satisfying \eqref{hyp2}, $a$ be a continuous function over $\bar\O$, and define
$$
M(x):=\int_{\O}m(x,y)\,dy,\quad \sigma:=\sup_{x\in\O}(a(x)-M(x)), \quad \Sigma:=\left\{y\in \bar\O|a(y)-M(y)=\sigma\right\}.
$$
Then, there exists a positive Radon measure  $d\mu_p% \in \mathbb{M}^+(\O)
$, such that, for any $\varphi$ in $C_c(\O)$, we have
$$
\int_{\O}\varphi(x)\left(\int_{\O}m(x,y)d\mu_p(y)\right) dx +\int_{\O}\varphi(x)(a(x)-M(x)+\lambda_p)d\mu_p(x)=0.
$$
%In addition, we have the following characterisation:
%\begin{itemize}
%\item Either there exists $varphi_p \in L^(\O)$, $\varphi_p>0$ such that $d\mu_p =\varphi_p(x)dx$
%\item Or there exists $ $
%\end{itemize}
%Let $\O$ be a bounded domain, $\m$ be defined by \eqref{bcl-eq-m} and assume that  $r\in C(\bar \O)$. Let  $M(x)$ and $\Sigma$ denote 
%$$ M(x):=\int_{\O}m(x,y)\,dy, \quad \Sigma:=\left\{y\in \bar\O|r(y)-M(y)=\sup_{x\in \O}(a(x)-M(x))\right\}.$$
%Then there exists a positive measure  $d\mu_p \in \M^+(\O)$, such that for any $\varphi\in C_c(\O)$ we have
%$$\int_{\O}\varphi(x)\left(\int_{\O}m(x,y)d\mu_p(y)\right) dx +\int_{\O}\varphi(x)(a(x)-M(x)+\lambda_p)d\mu_p(x)=0.$$
%Moreover, $d\mu_p$ can be  precisely characterised with re Either th there exists a positive function $\varphi_p \in L^1(\O)\cap C(\bar \O \setminus \Sigma)$ such that $\inf_{\O}\varphi_p>0$ and the measure $d\mu$ can be decomposed as follows $d\mu(x)= \varphi_p(x) dx + d\mu_s(x)$ where $d\mu_s(x)$ is a non negative singular measure with respect to the Lebesgue measure whose support lies in the set $\Sigma$.
In addition, we have the following dichotomy:
\begin{itemize}
\item either there exists $\varphi_p$ in $L^1(\O)$, $\varphi_p>0$, such that $d\mu_p =\varphi_p(x)dx$,
\item or there exists $g_p$ in $C(\bar\O)$, $g_p>0$, and $d\nu$ a positive singular measure with respect to the Lebesgue measure, whose support lies in the set $\Sigma$, such that
$$d\mu_p=\frac{g_p(x)}{\sigma-(a(x)-M(x))}dx +d\nu.$$
\end{itemize}
%there exists a positive function $g_p \in L^1(\O)\cap C(\bar \O \setminus \Sigma)$ such that $\inf_{\O}g_p>0$ and the measure $d\mu_p$ can be decomposed as follows $d\mu_p(x)= g_p(x) dx + d\mu_s(x)$ where $d\mu_s(x)$ is a non negative singular measure with respect to the Lebesgue measure whose support lies in the set $\Sigma$. 
\end{theorem} 
 
The measure $d\mu_p$ can be characterised more precisely and there exists a simple criterion guaranteeing its regularity.
  
\begin{proposition}[\cite{Coville2010,Coville2013,Berestycki2014,Coville2015}]\label{bcv-prop-phip}
Under the assumptions of the preceding theorem, $d\mu_p=\varphi_p(x)dx$ with   $\varphi_p\in C(\bar\O)$,  $\varphi_p>0$, if and only if $\ds{\lambda_p(\m+a)< -\sigma}$.  
%and $d\mu_p=\varphi_p(x)dx$. 
 \end{proposition}
 
We conclude by recalling a characterisation of $\lambda_p(\m +a)$ in the spirit of what is known for elliptic operators \cite{Berestycki2006, Berestycki2007, Berestycki2008}.  
%Motivated by the works \cite{Berestycki2006, Berestycki2007, Berestycki2008} on the generalised principal eigenvalue of an elliptic operators, $
%Let us introduce the following definition : 
\begin{definition}[\cite{Coville2013,Berestycki2014,Coville2015}]
Let $\O$ be a bounded domain, $\m$ be defined by \eqref{bcl-def-m} with a kernel satisfying \eqref{hyp2} and $a\in C(\bar \O)$. We define the following quantity:
\[
\lambda_p'(\m+a):=\inf\{\lambda\in\R \,|\, \exists \varphi\in C(\O)\cap L^\infty(\O),\ \varphi\ge 0,\ \opm{\varphi}+(a+\lambda)\varphi\geq 0\;\text{in}\; \O\}.
\]
%where $k(x):=\int_{\O}J(y-x)\,dy$ and $\langle \cdot, \cdot\rangle$ denotes the standard scalar product in $L^{2}(\O)$.
\end{definition}

As in the case of elliptic operators, the two quantities $\lambda_p$ and $\lambda_p'$ are  equal in our setting.% More precisely,

\begin{theorem}[\cite{Coville2013,Berestycki2014,Coville2015}]\label{bcl-thm-equal-lp}
 Let $\O$ be a bounded domain, $\m$ be defined by \eqref{bcl-def-m} with a kernel satisfying \eqref{hyp2} and $a\in C(\bar \O)$. Then 
$$
\lambda_p(\m +a)=\lambda_p'(\m +a).
$$
\end{theorem}

%It is worth to mention that this definition  has already been used for the study of \eqref{bcv-eq} in several papers \cite{Coville2008b,Coville2013,2013arXiv1305.7122C,Garcia-Melian2009a,Ignat2012}, but the relation between $\lambda_p,\lambda_p'$ and $\lambda_v$  had not been clarified.

%For elliptic operators, the analogues of these three quantities are equivalent on bounded domain \cite{Berestycki1994}. This is not necessarily the case for unbounded domains, where examples can be constructed \cite{Berestycki2006,Berestycki2007,Berestycki2010},  for which  $\lambda_1>\lambda_1'$. Since the operator, $\m+a$, shares many properties with elliptic operators, it is  suspected that the three quantities, $\lambda_p,\lambda_p'$ and $\lambda_v$, are not necessarily  equal.  However, for compactly supported kernel, $J$, we have:

\section{A priori estimates}\label{bcl-section-ae}
In this section, for a non-negative initial data $u_0$ in $L^1(\O)\cap L^{\infty}(\O)$, we establish some uniform in time \textit{a priori} estimates on the solution of \eqref{bcl-eq-intro}--\eqref{bcl-eq-intro-ci}. To do so, we start by proving a non-linear relative entropy identity satisfied by any solution of \eqref{bcl-eq-intro}. 
% We consider a parabolic  equation of the form
%  
%\begin{align}
%&\frac{\partial u}{\partial t}(t,x)=u(t,x)(a(x)-\Psi(x,u)(t))+\nabla\cdot(A(x)\nabla u(t,x) )\quad \text{ in } \quad \R_+\times\O, \label{bcl-eq-gen}\\
%&\frac{\partial u}{\partial n}(t,x)=0, \quad \text{ in } \quad \R_+\times\partial\O\label{bcl-eq-gen-bc}
%\end{align}
%where $\Psi(x,u)(t)$ denotes $\Psi(x,u)(t):=\int_{\O}K(x,y)|u|^p(t,y)\,dy.$
%Then for any solution of \eqref{bcl-eq-gen}--\eqref{bcl-eq-gen-bc} 

\begin{proposition}[general identity]\label{bcl-thm-gen-id}
Let $\O\subset \R^N$ be a bounded domain and assume that $a,k$ and $m$ satisfies \eqref{hyp1}--\eqref{hyp3}. Let $H$ be a smooth (at least $C^1$)  function. Let $\bar u$ be a $L^1(\O)\cap L^{\infty}(\O)$ positive stationary solution of \eqref{bcl-eq-intro}.  Let  $u\in C^1((0,+\infty),L^{\infty}(\O))$   be a  solution  of  \eqref{bcl-eq-intro}, then we have 
 
\begin{equation}
\frac{d \oph{H}{\bar u}{u}(t)}{dt}=-\D_{H,\bar u}(u)+\Gamma (t)\int_{\O} \bar u(x)H^{\prime}\left(\frac{u}{\bar u}(t,x)\right) u(t,x)\,dx
\end{equation}
where $\Gamma$, $\oph{H}{\bar u}{u}$, $\D_{H,\bar u}(u)$ are defined by
\begin{align*}
&\Gamma(t):=\int_{\O}k(y)(\bar u(y)-u(t,y))\,dy\\
&\oph{H}{\bar u}{u}(t):=\int_{\O}\bar u^2(x)H\left(\frac{u(t,x)}{\bar u(x) }\right)\, dx\\
&\D_{H,\bar u}(u)(t):=\iint_{\O\times \O} m(x,y)\bar u(x)\bar u(y)\left[H\left( \frac{u(t,x)}{\bar u(x)}\right)-H\left( \frac{u(t,y)}{\bar u(y)}\right) +H^{\prime}\left( \frac{u(t,x)}{\bar u(x)}\right) \left(\frac{u(t,x)}{\bar u(x)} -\frac{u(t,y)}{\bar u(y)}\right)\right] \, dxdy 
 \end{align*}
%where $(\vec a)^t $ denotes the transpose of a vector of  $\R^N$. 
\end{proposition}

\dem{Proof:}
From \eqref{bcl-eq-intro}, since the kernel $k$ satisfies condition \eqref{hyp3},  by defining $\Gamma(t):=\int_{\O}k(y)(\bar u(y)-u(t,y))\,dy$ we have for all $t>0$
\begin{equation*}
 \partial_t u(t,x)= \left(a(x) -\int_{\O}k(y)\bar u(y)\,dy\right)u(t,x)+\opm{u}(t,x)+\Gamma(t) u(t,x) \;\text{ for almost every }x\in\O.
\end{equation*}
Using that $\bar u$ is a positive stationary solution of \eqref{bcl-eq-intro}, for almost every  $x\in \O$, we have
$$
a(x)-\int_{\O}k(y)\bar u(y)\,dy =-\frac{1}{\bar u(x)}\opm{\bar u}(x),
$$
and we can rewrite the above equation  as follows
$$
\partial_t  u(t,x)=   \opm{u}(t,x) -\frac{u(t,x)}{\bar u(x)}\opm{\bar u}(x) +\Gamma(t) u(t,x) \;\text{ for almost every }\;x \in \O.
$$
Multiplying the above identity by $\bar u(x) H^{\prime}\left(\frac{u(t,x)}{\bar u(x)}\right)$ and integrating over $\O$, we find that
\begin{multline*}
\int_{\O} \bar u(x)H^{\prime}\left(\frac{u(t,x)}{\bar u(x)}\right) \partial_t u(t,x)\, dx=\int_{\O} \bar u(x)H^{\prime}\left(\frac{u(t,x)}{\bar u(x)}\right) \Gamma(t) u(t,x)\, dx\\+\int_{\O} H^{\prime}\left(\frac{u(t,x)}{\bar u(x)}\right)[\bar u(x)\opm{u}(t,x)-u(t,x)\opm{\bar u}(x)]\, dx .
\end{multline*}
By rearranging the terms, we get 
\begin{multline*}
\int_{\O} \bar u(x)H^{\prime}\left(\frac{u(t,x)}{\bar u(x)}\right) \partial_t u(t,x)\, dx=\int_{\O} \bar u(x)H^{\prime}\left(\frac{u(x)}{\bar u(t,x)}\right) \Gamma(t) u(t,x)\, dx\\ + \iint_{\O\times\O} m(x,y)\bar u(x)\bar u(y)H^{\prime}\left( \frac{u(t,x)}{\bar u(x)}\right) \left[\frac{u(t,y)}{\bar u(y)}-\frac{u(t,x)}{\bar u(x)} \right]\,dxdy,
\end{multline*}
and, due to the symmetry of $m$, we straightforwardly see that 
$$
\iint_{\O\times\O} m(x,y)\bar u(x)\bar u(y)\left[H\left(\frac{u(t,x)}{\bar u(x)}\right)-H\left(\frac{u(t,y)}{\bar u(y)}\right) \right] \, dxdy =0.
$$
Hence, by combining the above equalities, we reach  
$$
\frac{d }{dt}\oph{H}{\bar u}{u}(t)=\Gamma(t)\int_{\O} \bar u(x)H^{\prime}\left(\frac{u(t,x)}{\bar u(x)}\right)u(t,x)\, dx-\D_{H,\bar u}(u).
$$
\fdem
%\begin{remark}\label{bcl-rem-id}
%We want to stress that if we replace $\bar u$ by any positive function $\tilde u$ satisfying
%\begin{align*}
%&\nabla\cdot(A(x)\nabla \tilde u(x) )=-\tilde u(x)\left(a(x)-\tilde \Psi(x,\tilde u)(t)\right)\quad \text{ in } \quad \O, \\
%&\frac{\partial \tilde u}{\partial n}(x)=0, \quad \text{ in } \quad \partial\O
%\end{align*}
%it will affect the equality in Theorem \ref{bcl-thm-gen-id} only through the term $\Gamma$  which will be transform into
% $$\Gamma(t,x) =\tilde \Psi(x,\tilde u(x))-\Psi(x,u(t,x)).$$
%\end{remark}
%\begin{remark}
%Under the  extra assumption $\frac{u}{\bar u} \in L^{\infty}(\O)$, we remark that the formulas will holds as well if we consider  homogeneous Dirichlet boundary conditions instead of Neumann boundary conditions. It is worth noticing that this extra condition is always satisfied in the Neumann case since for all positive stationary solution with homogeneous Neumann Boundary condition, we can show that  $\inf_{\bar \O}\bar u>0$.
%\end{remark}

\begin{remark}
When $k\equiv 0$, equation \eqref{bcl-eq-intro} is linear and relative entropy formulas  are well known in this case, see \cite{Michel2005}. 
\end{remark}
\begin{remark}\label{bcl-rema-gen-id}
When $H$ is non decreasing and $\bar u$ is only assumed to be a stationary super-solution of \eqref{bcl-eq-intro},  from the above proof we clearly see that, 
$$\frac{d }{dt}\oph{H}{\bar u}{u}(t)\le -\D_{H,\bar u}(u)(t)+\Gamma(t)\int_{\O} \bar u(x)H^{\prime}\left(\frac{u}{\bar u}(t,x)\right) u(t,x)\,dx.
$$
Similarly, if $\bar u$ is a positive stationary sub-solution of \eqref{bcl-eq-intro}, we have   
$$\frac{d }{dt}\oph{H}{\bar u}{u}(t)\ge -\D_{H,\bar u}(u)(t)+\Gamma(t)\int_{\O} \bar u(x)H^{\prime}\left(\frac{u}{\bar u}(t,x)\right) u(t,x)\,dx.
$$
\end{remark}

Equipped with this general relative entropy identity,  we may derive  some  useful differential inequalities.

\begin{proposition}\label{bcl-lem-liap}
Let $\O\subset \R^N$ be a bounded domain and assume that $a,k$ and $m$ satisfies \eqref{hyp1}--\eqref{hyp3}. Let $q\ge 1$ and $H_q$ be the smooth convex function $H_q(s):\,s\mapsto s^q$. Let $u,\bar u$ be two positive  solutions of \eqref{bcl-eq-intro} as in Proposition \ref{bcl-thm-gen-id}. 
%Assume further that  $\bar u$ is a stationary solution of \eqref{bcl-eq-intro}. 
Then the functional  $\f(u):=\log\left(\frac{\oph{q}{\bar u}{u}}{\left(\oph{1}{\bar u}{u}\right)^{q}}\right)$ satisfies: 
\begin{equation} 
\frac{d}{dt}\f(u)(t)=-\frac{1}{\oph{q}{\bar u}{u}(t)} \D_{q,\bar u}(u)(t)\le 0.\label{bcl-eq-liap}
\end{equation}

%Then, we have
%\begin{equation}
%\frac{d \oph{q}{\bar u}{u}(t)}{dt}=-\D_{q,\bar u}(u) + q\Gamma(t)\oph{q}{\bar u}{u}(t). \label{bcl-eq-estim1}
%\end{equation}
%where $\ds{\oph{q}{\bar u}{u}(t):=\int_{\O}\bar u^2(x)\left(\frac{u(t,x)}{\bar u(x)} \right)^q\,dx}$. 
Moreover, we have 
\begin{equation}
\D_{q,\bar u}(u)(t)=\frac{q}{2}\iint_{\O\times\O}m(x,y)\bar u(x)\bar u(y)\left[\frac{u(t,x)}{\bar u(x)}-\frac{u(t,y)}{\bar u(y)}\right]\left[\left(\frac{u(t,x)}{\bar u(x)}\right)^{q-1}-\left(\frac{u(t,y)}{\bar u(y)}\right)^{q-1}\right]\,dxdy. \label{bcl-estim}
\end{equation}

%Furthermore, the functional  $\f(u):=\log\left(\frac{\oph{q}{\bar u}{u}(t)}{\left(\oph{1}{\bar u}{u}(t)\right)^{q}}\right)$ satisfies: 
%\begin{equation} \frac{d}{dt}\f(u)=-\frac{1}{\oph{q}{\bar u}{u}(t)} \D_{q,\bar u}(u)\le 0.\label{bcl-eq-liap}\end{equation}
%where $\oph{1}{\bar u}(u):=\int_{\O}\bar u u$.
\end{proposition}
\medskip

\begin{remark}
We observe that in the case of $H(s)=s^2$, $\oph{2}{\bar u}{u}(t)=\nlto{u(t,\cdot)}^2$ and  we get a Lyapunov functional involving the $L^2$ norm of $u$ instead of a weighted $L^q$ norm of $u$. Indeed, we have
$$
\frac{d }{d t}\left(\log\left(\frac{\nlto{u(t,\cdot)}^2}{\left(\oph{1}{\bar u}{u} (t)\right)^2}\right)\right)=-\frac{1}{\nlto{u(t,\cdot)}^2}\iint_{\O\times \O}m(x,y)\bar u(x)\bar u(y)\left(\frac{u(t,x)}{\bar u(x)}-\frac{u(t,y)}{\bar u(y)}\right)^2\,dxdy.
$$
%In this particular situation,  the assumption on $ \bar u$ can be relaxed and   by  assuming that $\bar u$ is either  super or a sub-solution of \eqref{bcl-eq-stat} we get differential inequalities instead of equalities. For example, if $\bar u$ is a sub-solution, we get
%$$ \frac{d }{d t}\left(\log\left(\frac{\nlto{u}^2}{\left(\oph{1}{\bar u}{u} \right)^2}\right)\right)\le -\frac{1}{\nlto{u}^2}\iint_{\O\times \O}m(x,y)\bar u(x)\bar u(y) \left( \frac{u(t,x)} {\bar u(x)}-\frac{u(t,y)} {\bar u(y)}\right)^2\, dydx.$$
\end{remark}

\dem{Proof of Proposition \ref{bcl-lem-liap}:}
Observe that, for $H(s):=s^q$, we have, by Proposition \ref{bcl-thm-gen-id},
$$ 
\frac{d}{dt}\oph{q}{\bar u}{u}(t)=-\D_{q,\bar u}(u)(t)+q\,\Gamma(t)\int_{\O}\bar u (x)\left(\frac{u(t,x)}{\bar u(x)}\right)^{q-1}u(t,x)\,dx.
$$
Therefore, by definition of $\oph{q}{\bar u}{u}$ we get from the above equality
\begin{equation}\label{bcl-eq-estim1}
\frac{d }{dt}\oph{q}{\bar u}{u}(t)=-\D_{q,\bar u}(u)(t)+q\,\Gamma(t)\oph{q}{\bar u}{u}(t).
\end{equation}
Now by taking $q=1$ in \eqref{bcl-eq-estim1}, we obtain
\[
\frac{d }{dt}\oph{1}{\bar u}{u}(t)=-\D_{1,\bar u}(u)(t)+\Gamma(t)\oph{1}{\bar u}{u}(t).
\]
A quick computation shows that $\ds{\D_{1,\bar u}(u) }=0$ and therefore 
\begin{equation}\label{bcl-eq-estim2}
\frac{d }{dt}\oph{1}{\bar u}{u}(t)=\Gamma(t)\oph{1}{\bar u}{u}(t).
\end{equation}
%The identity   \eqref{bcl-eq-estim1} is a straightforward consequence of Lemma \ref{bcl-thm-gen-id}. Indeed, for $H(s):=s^q$, by the Theorem \ref{bcl-thm-gen-id} we have:
%$$ 
% \frac{d \oph{q}{\bar u}{u}}{dt}(t)=-\D_{q,\bar u}(u)+\int_{\O} \bar u (x)H^{\prime}\left(\frac{u(t,x)}{\bar u(x)}\right) \Gamma (x)u(x)\,dx
%$$
% where $\Gamma$, $\D$ are the following quantity:
%\begin{align*}
%&\Gamma(u(t)):=\Psi(\bar u) -\Psi(u)\\
%&\D(u):=\int_{\O} H^{\prime\prime}\left( \frac{u(t,x)}{\bar u(x)}\right) \bar u^2(x) \left(\nabla\left( \frac{u(t,x)} {\bar u(x)}\right)\right)^tA(x)\nabla\left( \frac{u(t,x)} {\bar u(x)}\right) \, dx
%\end{align*}
%By observing that $\bar u(x) u(x)H^{\prime}\left(\frac{u(t,x)}{\bar u(x)}\right)=q\oph{q}{\bar u}{u}(t)$ and that $\Gamma$ is independent of $x$, we see that 
% $$\frac{d \oph{q}{\bar u}{u}(t)}{dt}=-\D(u)+q\Gamma \oph{q}{\bar u}{u}(t),$$
%and the formula \eqref{bcl-eq-estim1} holds.
%To obtain \eqref{bcl-eq-liap},  we observe that by taking $q=1$ in the  formula \eqref{bcl-eq-estim1} we get    
%$$ \frac{d \oph{1}{\bar u}{u}(t)}{dt}=\Gamma\oph{1}{\bar u}{u}(t).$$
Since for all $q\ge 1$, $\oph{q}{\bar u}{u}(t)>0$ for all times, we have
\begin{align} 
&\frac{d}{dt}\log(\oph{1}{\bar u}{u}(t))=\Gamma(t),\label{bcl-eq-cdif1}\\
&\frac{d}{dt}\log(\oph{q}{\bar u}{u}(t))=-\frac{1}{\oph{q}{\bar u}{u}(t)}\D_{q,\bar u}(u)(t)+q\,\Gamma(t).\label{bcl-eq-cdif2}
\end{align}
By combining \eqref{bcl-eq-cdif1} and \eqref{bcl-eq-cdif2}, we end up with 
$$
\frac{d }{d t}\left(\log\left(\frac{\oph{q}{\bar u}{u}(t)}{\left(\oph{1}{\bar u}{u}(t)\right)^{q}}\right)\right)=-\frac{1}{\oph{q}{\bar u}{u}(t)}\D_{q,\bar u}(u)(t).
$$
Equality \eqref{bcl-estim} then follows straightforwardly from direct computations, by using symmetry and an obvious change of variables. 

%\begin{multline*}
%\D_{q,\bar u}(u)=\frac{1}{2}\iint_{\O\times\O}m(x,y)\bar u(x)\bar u(y)\left[q\left(\frac{u(t,x)}{\bar u(x)}\right)^{q-1}\left(\frac{u(t,x)}{\bar u(x)}-\frac{u(t,y)}{\bar u(y)}\right)\\+q\left(\frac{u(t,y)}{\bar u(y)}\right)^{q-1}\left(\frac{u(t,y)}{\bar u(y)}-\frac{u(t,x)}{\bar u(x)}\right)\right ].
%\end{multline*}
\fdem 

%\begin{remark}
%We observe that in the case of $H(s)=s^2$, $\oph{2}{\bar u}{u}=\nlto{u}^2$ and  we get a Lyapunov functional involving the $L^2$ norm of $u$ instead of  a weighted $L^q$ norm of $u$. 
%Indeed, we have $$ \frac{d }{d t}\left(\log\left(\frac{\nlto{u}^2}{\left(\oph{1}{\bar u}{u} \right)^2}\right)\right)=-\frac{1}{\nlto{u}^2}\iint_{\O\times \O}m(x,y)\bar u(x)\bar u(y) \left( \frac{u(t,x)} {\bar u(x)}-\frac{u(t,y)} {\bar u(y)}\right)^2\, dydx.$$
%%In this particular situation,  the assumption on $ \bar u$ can be relaxed and   by  assuming that $\bar u$ is either  super or a sub-solution of \eqref{bcl-eq-stat} we get differential inequalities instead of equalities. For example, if $\bar u$ is a sub-solution, we get
%%$$ \frac{d }{d t}\left(\log\left(\frac{\nlto{u}^2}{\left(\oph{1}{\bar u}{u} \right)^2}\right)\right)\le -\frac{1}{\nlto{u}^2}\iint_{\O\times \O}m(x,y)\bar u(x)\bar u(y) \left( \frac{u(t,x)} {\bar u(x)}-\frac{u(t,y)} {\bar u(y)}\right)^2\, dydx.$$
%\end{remark}

From these differential inequalities, we obtain  uniform in time \textit{a priori} bounds of the $L^1$ norm of a solution of \eqref{bcl-eq-intro} -- \eqref{bcl-eq-intro-ci}. Namely, we show  
%\begin{lemma}\label{bcl-lem-esti}
% Let $u\in C^1((0,+\infty), L^{2}(\O)\cap L^1(\O))$ be a non-negative solution of the Cauchy problem \eqref{bcl-eq} with initial data $u_0\in L^2(\O)\cap L^1(\O)$, $u_0\ge 0$.   Then  for all $1\le q\le 2$ there exists  two positive constants $C_q(q,u_0)>c_q(q,u_0)>0$ such that for all $t \ge 0$ 
% $$c_q\le \nlp{u}{q}{\O}\le C_q. $$
%\end{lemma}
\begin{lemma}\label{bcl-lem-esti}
Let $\O\subset \R^N$ be a bounded domain and assume that $a,k$ and $m$ satisfies \eqref{hyp1}--\eqref{hyp3}. Let $u\in C^1((0,+\infty), L^{1}(\O)\cap L^{\infty}(\O))$ be a non-negative solution of the Cauchy problem \eqref{bcl-eq-intro}--\eqref{bcl-eq-intro-ci} with initial data $u_0\in L^{\infty}(\O)\cap L^1(\O)$, $u_0\ge 0$.   Then  there exists  two positive constants $C_1(u_0)>c_1(u_0)>0$ such that for all $t \ge 0$ 
 $$c_1\le \nlp{u(t,\cdot)}{1}{\O}\le C_1. $$
\end{lemma}
\dem{Proof:}
%First, we observe that it is enough to obtain the bounds for $q=1$ and $q=2$, the estimates for $\in (1,2)$ is then  obtained by interpolation.

%Let us first show that for all $q\ge 1$ then there exists $C_q(q,u(x,1))$ so that for all $t \ge 1$ 
% \begin{equation}
% \nlp{u}{q}{\O}\le C_q. \label{bcl-eq-blind-esti1}
% \end{equation}
%Let us start with  the  \textit{ a priori}  bounds  for $\nlp{u}{1}{\O}$.  
First, let us observe that large (respectively small) constants are super-solutions (respectively sub-solutions) of \eqref{bcl-eq-intro}. Indeed, for $\ds{C\ge \sup_{x\in \O}\left(\frac{a(x)}{\int_{\O}k(z)\,dz} \right)}$, we have 
$$
\opm{C}+C\left(a(x)-C\int_{\O}k(z)\,dz\right)= C(a(x)-\sup_{x\in \O}a(x))\le 0.
$$
Similarly, for $\ds{C\le \inf_{x\in \O}\left(\frac{a(x)}{\int_{\O}k(z)\,dz} \right)}$, we get
$$
\opm{C}+C \left(a(x)-C\int_{\O}k(z)\,dz\right)= C(a(x)-\inf_{x\in \O}a(x))\ge 0.
$$

Therefore, from Proposition \ref{bcl-thm-gen-id} and Remark \ref{bcl-rema-gen-id}, by choosing $\bar u$ a large, respectively a small constant, and considering the convex function $H(s):s\mapsto s$, we get  
\begin{align*}
\frac{d}{d t }(\nlp{u(t,\cdot)}{1}{\O}) \le C\left(\int_{\O}k(z)\,dz - \int_{\O}k(z)u(t,z)\,dz\right)\nlp{u(t,\cdot)}{1}{\O}\quad \text{ for } \quad C\ge  \sup_{x\in \O}\left(\frac{a(x)}{\int_{\O}k(z)\,dz} \right),
\\
\frac{d}{d t }(\nlp{u(t,\cdot)}{1}{\O})\ge c\left(\int_{\O}k(z)\,dz - \int_{\O}k(z)u(t,z)\,dz\right)\nlp{u(t,\cdot)}{1}{\O}\quad \text{ for } \quad c\le  \inf_{x\in \O}\left(\frac{a(x)}{\int_{\O}k(z)\,dz} \right).
\end{align*}

Now, since $k$ satisfies \eqref{hyp3}, we have  $c_0\le k(y)\le C_0$ for all $y\in \O$ and from the above differential inequalities we get 
\begin{align*}
\frac{d}{d t }(\nlp{u(t,\cdot)}{1}{\O}) \le C\left(\int_{\O}k(z)\,dz - c_0\nlp{u(t,\cdot)}{1}{\O}\right)\nlp{u(t,\cdot)}{1}{\O}\quad \text{ for } \quad C\ge  \sup_{x\in \O}\left(\frac{a(x)}{\int_{\O}k(z)\,dz} \right),
\\
\frac{d}{d t }(\nlp{u(t,\cdot)}{1}{\O})\ge c\left(\int_{\O}k(z)\,dz - C_0\nlp{u(t,\cdot)}{1}{\O}\right)\nlp{u(t,\cdot)}{1}{\O}\quad \text{ for } \quad c\le  \inf_{x\in \O}\left(\frac{a(x)}{\int_{\O}k(z)\,dz} \right).
\end{align*}
From the logistic character of these two differential inequalities and since $\nlp{u_0}{1}{\O}>0$, we deduce that for all $t\ge 0$  
$$c_1:=\inf\left\{\nlp{u_0}{1}{\O}; \frac{\int_{\O}k(z)\,dz}{C_0}\right\} \le u(t,x)\le C_1:=  \sup\left\{\nlp{u_0}{1}{\O}; \frac{\int_{\O}k(z)\,dz}{c_0}\right\}.$$

\fdem

\begin{remark}
Observe that the above proof holds as well for $k$ bounded above and below by positive constants. As a consequence, such a uniform $L^1$ estimate can be also obtained in a more general situation where the kernel $k$ is not necessarily independent of the trait $x$.    
\end{remark}

\section{Proofs}\label{bcl-section-proof}
We are now in a position to prove Theorems \ref{bcl-thm1} and \ref{bcl-thm2}. Let us start with the construction of a stationary measure.

\subsection{Construction of a Stationary state}
Consider the stationary problem \eqref{bcl-eq-red}, then in order to construct stationary state in the space of Radon measures, we have to find $d\mu$  solution of the following weak formulation 
\begin{equation} \label{bcl-eq-redw-proof}
\forall\varphi \in C_c(\O),\ \int_{\O}\left(\opm{\varphi}(x)+a(x)\varphi(x)\right)\,d\mu(x)=\int_{\O}\varphi(x)d\mu(x)\int_{\O}k(x)d\mu(x).  
\end{equation}
 
%This situation appears for example  when $\ds{\lambda_p(\M+r)<-\sigma}$, (Theorem \ref{bcl-thm-critex}). In this situation, there exists a positive continuous function $\varphi_p$ associated with $\lambda_p(\M+r)$. To prove Theorem \ref{bcl-thm}, we start by constructing a stationary solution.
Owing to Theorem \ref{bcl-thm-mu}, let us consider a positive measure $d\mu_p$ associated to $\ds{\lambda_p(\m+a)}$, which we normalise in order to have $\int_{\O}d\mu_p=1$.% We claim that  
%In this situation, there exists a positive continuous function $\varphi_p$ associated with $\lambda_p(\M+r)$. To prove Theorem \ref{bcl-thm}, we start by constructing a stationary solution.

%\subsection*{Stationary solution}
 %In this situation, we claim 
\begin{claim}\label{bcl-cla-steady}
There exists a unique $\theta>0$ such that $\theta\,d\mu_p$ is a  positive stationary solution  of \eqref{bcl-eq-redw-proof}. 
\end{claim}
\dem{Proof:}
Let $\theta$ be defined by
$\theta:=\left(\frac{-\lambda_p(\m +a)}{\int_{\O}k(y)d\mu_p(y)\,dy}\right)$. For any $\varphi\in C_c(\O)$, $\theta\,d\mu_p$ satisfies 
\begin{align*}
\int_{\O}\left(\opm{\varphi}(x)+a(x)\varphi(x) \right) \theta\,d\mu_p(x)&=-(\lambda_p(\m+a)) \int_{\O}\varphi(x)\theta\,d\mu_p(x)\\
&=-\frac{(\lambda_p(\m+a))}{\int_{\O}k(x)\,d\mu_p(x)}\int_{\O}k(x)\,d\mu_p(x) \int_{\O}\varphi(x)\theta\,d\mu_p(x)\\
&=\int_{\O}k(x)\theta\,d\mu_p(x) \int_{\O}\varphi(x)\theta\,d\mu_p(x).
\end{align*} 
 
Thus, $\theta\,d\mu_p$ is a stationary solution of \eqref{bcl-eq-redw-proof}.

To conclude, it remains to show that  $-\lambda_p(\m+a)>0$. This is the case, since we have $\lambda_p(\m+a)=\lambda_p'(\m+a)$ by Theorem \ref{bcl-thm-equal-lp}, and by taking $(-\inf_{x\in \O} a(x),1)$ as test function, we can easily check that $$\lambda_p'(\m+a)\le -\inf_{x\in \O} a(x)<0.$$
\fdem

\begin{remark}
From the above computation, we clearly see  that the uniqueness of the stationary state follows from the uniqueness of the  measure associated with $\lambda_p$.  
\end{remark}
 
\subsection{Long time behaviour}
Let us now prove Theorem \ref{bcl-thm2}. We assume that $d\mu$ the positive measure constructed above is regular and bounded, i.e. $d\mu(x)=\bar u(x)\,dx$ with $\bar u \in L^{1}(\O)\cap L^{\infty}(\O)$. Since $d\mu$ is associated with the principal eigenvalue $\lambda_p(\m+a)$, then $\bar u=\theta \varphi_p$ with $\varphi_p\in L^1(\O)\cap L^{\infty}(\O)$ and $\theta$ defined in the proof of Claim \ref{bcl-cla-steady}. From the regularity of $\varphi_p$, we can see that $\bar u$ is a strong solution of \eqref{bcl-eq-red}. Now, knowing that a positive continuous stationary solution of \eqref{bcl-eq-red} exists, we can derive further \textit{a priori } estimates on the solution $u$ of \eqref{bcl-eq-intro}--\eqref{bcl-eq-intro-ci}.

\begin{lemma}\label{bcl-lem-esti2}
Let $\O\subset \R^N$ be a bounded domain and assume that $a$, $k$ and $m$ satisfy \eqref{hyp1}--\eqref{hyp3}. Let $u\in C^1((0,+\infty),L^{1}(\O)\cap L^{\infty}(\O))$ be a non-negative solution of the Cauchy problem \eqref{bcl-eq-intro}--\eqref{bcl-eq-intro-ci} with initial datum $u_0\in L^{\infty}(\O)\cap L^1(\O)$, $u_0\ge 0$.  Then there exist two constants $C_2>c_2>0$ depending on $u_0$ such that for all $t>0$,
$$
c_2\le\nlp{u(t,\cdot)}{2}{\O}< C_2.
$$
\end{lemma}
\dem{Proof:}
The uniform lower bound is rather easy to obtain and follows directly from the  H{\"o}lder's inequality and the estimates in Lemma \ref{bcl-lem-esti}. Indeed, since $\O$ is bounded we have 
$$
\sqrt{c_1}\le\nlp{u(t,\cdot)}{1}{\O}\le\sqrt{|\O|}\nlp{u(t,\cdot)}{2}{\O},
$$
where $|\O|$ denotes the Lebesgue measure of the set $\O$.

On the other hand, we get an uniform upper bound  as a straightforward application of  Proposition \ref{bcl-lem-liap}. Namely, since $\bar u$ is a positive $L^1(\O)\cap L^{\infty}(\O)$ stationary solution of \eqref{bcl-eq-intro}, by Proposition \ref{bcl-lem-liap}, the functional $F(t):=2\log\left(\frac{\nlp{u(t,\cdot)}{2}{\O}}{\oph{1}{\bar u}{u}(t)}\right) $ is a decreasing function of $t$ and therefore, for all $t\ge 0$,
$$
\nlp{u(t,\cdot)}{2}{\O}\le \oph{1}{\bar u}{u}(t)\left(\frac{\nlp{u_0}{2}{\O}}{\oph{1}{\bar u}{u_0}}\right)\le C(u_0)\|\bar u\|_{\infty} \nlp{u(t,\cdot)}{1}{\O}. $$
Hence, for all $t\ge 0$,
$$
\nlp{u(t,\cdot)}{2}{\O}\le C_1C(u_0)\|\bar u\|_{\infty}.
$$
\fdem

In order to prove that $u$ converges to a stationary solution, we introduce the following decomposition of $u$.
%Since $u(t,x)\in C^{1}((0,+\infty), L^{1}(\O)\cap L^{\infty})$ be a positive solution of \eqref{bcl-eq} and let   $\bar u$ be  the stationary solution  constructed in the previous Claim \ref {bcl-cla-steady} , i.e $\bar u:=\theta\varphi_p$. 
Since for all $t>0$, $u$ and $\bar u$ belong to $ L^{1}(\O)\cap L^{\infty}(\O)$, they belong to $L^2(\O)$ and we can write $u$ as follows:
$$
u(t,x):=\lambda(t) \bar u(x) +h(t,x)
$$
with $h$ such that $\int_{\O}\varphi_p(x)h(t,x)\,dx=0$ for all $t>0$.

%Equipped with this bounded, we are now in position to show that $\lambda(t)\to 1$ and $h\to 0$ and finish the proof of the Theorem \ref{bcl-thm2}.

%We claim 
\begin{claim}
$\lambda(t)\to 1$ and $\nlp{h(t,\cdot)}{2}{\O}\to 0$ as $t\to+\infty$.
\end{claim}
\dem{Proof:}
For convenience, we introduce the following notation $\la\varphi,\psi\ra:=\int_{\O}\varphi(x)\psi(x)\,dx$ to denote the standard scalar product of two function of $L^2(\O)$.

We start by deriving some useful bounds on $\lambda$ and $h$. From the decomposition, we have 
$$
\la \bar u, u(t,\cdot)\ra = \lambda(t) \theta^2=\theta \int_{\O}u(t,x)\varphi_p(x)\,dx. 
$$
Therefore, since $\varphi_p$ is positive and bounded in $\bar \O$, we have from Lemma \ref{bcl-lem-esti}
$$
c_1(u_0)\inf_{x\in \O}\varphi_p(x)\le \int_{\O}u(t,x)\varphi_p(x)\,dx \le C_1(u_0)\|\varphi_p\|_{\infty},
$$
and 
\begin{equation}\label{bcl-eq-esti-lpt}
 \frac{c_1(u_0)\inf_{x\in \O}\varphi_p(x)}{\theta}\le \lambda(t) \le \frac{C_1(u_0)\|\varphi_p\|_{\infty}}{\theta}.
\end{equation}

From Lemma \ref{bcl-lem-esti2}, we obviously derive an upper bound for $\nlp{h(t,\cdot)}{2}{\O}$. Indeed, by construction
\begin{equation}\label{bcl-eq-esti-h}
\nlp{h(t,\cdot)}{2}{\O}^2\le \nlp{u(t,\cdot)}{2}{\O}^2\le C_2.
\end{equation}

Substituting to $u$ its decomposition in the equation \eqref{bcl-eq-intro}, we get 
\begin{align}
\lambda'(t)\bar u(x) +\frac{\partial }{\partial t} h(t,x)
&=\left(a(x) -\int_{\O}k(z)u(t,z)\,dz\right)\left(\lambda(t)\bar u(x) +h(t,x)\right) +\lambda(t)\opm{\bar u}(x)+\opm{h}(t,x). \label{bcl-eq-red-decomp1}
\end{align} 
Multiplying the above equation by $h$ and integrating it over $\O$, we get after obvious computations
$$
\left\langle\frac{\partial}{\partial t}h(t,\cdot),h(t,\cdot)\right\rangle= \left\langle\left(a(x) -\int_{\O}k(z)u(t,z)\,dz\right)h(t,\cdot)+\opm{h}(t,\cdot),h(t,\cdot)\right\rangle,
$$
where we used the definition of $\bar u$ and $\la \bar u,h(t,\cdot)\ra=0$.

Since $\oph{2}{\bar u}{h}(t):=\nl{h(t,\cdot)}{\O}^2$, we get 
$$
\left\langle\frac{\partial h}{\partial t}(t,\cdot),h(t,\cdot)\right\rangle=\frac{1}{2}\frac{d}{dt}\oph{2}{\bar u}{h}(t)=\left\langle\left(a(x) -\int_{\O}k(z)u(t,z)\,dz\right)h(t,\cdot)+\opm{h}(t,\cdot),h(t,\cdot)\right\rangle.
$$

By following the computation developed for the proof of Proposition \ref{bcl-thm-gen-id} with $H(s)=s^2$, we see that 
\begin{equation}\label{bcl-eq-h}
\frac{d}{dt}\oph{2}{\bar u}{h}(t)=-\D_{2,\bar u}(h)(t)+\Gamma(t)\oph{2}{\bar u}{h}(t).
\end{equation}
with $\ds{\Gamma(t):=\left(\int_{\O}k(z)\bar u(z)\,dz -\int_{\O}k(z)u(t,z)\,dz\right)=\left(-\lambda_p -\int_{\O}k(z)u(t,z)\,dz\right)}$.

By construction, $\oph{2}{\bar u}{h}(t)\ge 0$ for all $t\ge 0$. So either $\oph{2}{\bar u}{h}(t)>0$ for all times $t$ or  there exists $t_0\in \R$ such that $\oph{2}{\bar u}{h}(t_0)=0$. In the latter case, we have $u(t_0,x)=\lambda(t_0)\theta \varphi_p(x)$ for almost every $x\in \O$. Let  $w(t,x):=\gamma(t)\theta \varphi_p$ with $\gamma(t)$ satisfying the ODE
\begin{align}
&\frac{d}{dt}\gamma(t)=-\lambda_p \gamma(t)(1-\gamma(t))\\
&\gamma(t_0)=\lambda(t_0).
\end{align}
By construction, $\gamma(t)\to 1$ as $t\to+\infty$ and we can check that $w$ is a solution of \eqref{bcl-eq-intro} for all $t\ge t_0$. Thus, since $w(t_0,\cdot)=u(t_0,\cdot)$ by uniqueness of the solution of the Cauchy problem \eqref{bcl-eq-intro}, we have $u(t,\cdot)\equiv w(t,\cdot)$ for all $t\ge t_0$ and therefore for all $t\ge t_0$, $h(t,\cdot)\equiv 0$ and $\lambda(t)=\gamma(t)$.

In the other situation, $\oph{2}{\bar u}{h}(t)>0$ for all $t$ and we claim the following.
 
\begin{claim} \label{bcl-cla-energy}
$\oph{2}{\bar u}{h}(t)\to 0$ as $t\to +\infty$.
\end{claim}

Assume the Claim holds then we can conclude the proof by arguing as follows. From the decomposition $u(t,x)=\lambda(t)\bar u(x) +h(t,x)$, we can express  the function  $\oph{1}{\bar u}{u}(t)$ by $\oph{1}{\bar u}{u}(t)=<\bar u,u(t,\cdot)>=\lambda(t)\la \bar u ,\bar u\ra$. Using Proposition \ref{bcl-lem-liap}, we deduce that 
\begin{equation}
\lambda'(t)=-\lambda_p(1-\lambda(t))\lambda(t)-\left(\int_{\O}k(z)h(t,z)\,dz\right)\lambda(t). \label{bcl-eq-lambda}
\end{equation}

Now by using  $\nlto{h(t,\cdot)}^2=\oph{2}{\bar u}{h}(t)\to 0$ as $t\to+\infty$, we deduce that
$$
\lambda'(t)=-\lambda_p(1-\lambda(t))\lambda(t)+ \lambda(t)\,o(1),
$$
with $o(1)\le C\nlp{h(t,\cdot)}{2}{\O}$.

Therefore, by a elementary  analysis of the ODE, we deduce that $\lambda(t)\to 1$ as $t\to+\infty$. 
%which prove that  $u$ converges  to $\bar u$ almost everywhere.
\fdem

\dem{Proof of Claim \ref{bcl-cla-energy}:} 
Since $\oph{2}{\bar u}{h}(t)>0$ for all $t$,  from \eqref{bcl-eq-h} and by following the proof of Proposition \ref{bcl-lem-liap} we see that 
\begin{equation}
\frac{d}{dt}\log\left[\frac{\oph{2}{\bar u}{h}(t)}{\left(\oph{1}{\bar u}{u}(t)\right)^2}\right]=  -\frac{\D_{2,\bar u}(h)(t)}{\oph{2}{\bar u}{h}(t)} \le 0\label{bcl-cla-dF}
\end{equation} 
Thus $ F(t):= \log\left[\frac{\oph{2}{\bar u}{h}(t)}{\left(\oph{1}{\bar u}{u}(t)\right)^2}\right]$ is a non increasing smooth function.
Thanks the monotonicity of $F$, to prove the Claim, it is sufficient to exhibit a sequence $(t_n)_{n\in\N}$ such that $t_n\to +\infty$ and  $ \oph{2}{\bar u}{h}(t_n)\to 0$. 

%Indeed, assume such sequence exists and let  $(s_k)_{k\in \N}$ be a sequence going to $+\infty$.  Then  there exists $k_0$ and a subsequence $(t_{n_k})_{k\in \N}$ of  $(t_n)_{n\in\N}$ so that  for all $k\ge k_0$,  we have $s_k \ge t_{n_k}$. Therefore from the monotonicity of $\tilde F$ we have for all $k\ge k_0$
%$$\log\left[\frac{\oph{2}{\bar u}{h}(s_k)}{\left(\oph{1}{\bar u}{u}(s_k)\right)^2}\right]\le \log\left[\frac{\oph{2}{\bar u}{h}(t_{n_k})}{\left(\oph{1}{\bar u}{u}(t_{n_k})\right)^2}\right].$$
%By letting $k$ to infinity in the above inequality, we deduce that 
%$$\lim_{k\to \infty}\log\left[\frac{\oph{2}{\bar u}{h}(s_k)}{\left(\oph{1}{\bar u}{u}(s_k)\right)^2}\right]=-\infty,$$
%which implies that $\oph{2}{\bar u}{h}(s_k)\to 0, $ since by Lemma \ref{bcl-lem-esti}  $(\oph{1}{\bar u}{u}(t_k))_{k\in \N}$ is uniformly bounded.
%The sequence $(s_k)_{k \in \N}$ being chosen arbitrarily this implies that $\oph{2}{\bar u}{h}(t) \to 0 $ as $t\to +\infty$.  
%\medskip

To exhibit such sequence, it is sufficient to prove that $\inf_{t\in \R_+}\oph{2}{\bar u}{h}(t)=0$. By contradiction, let us assume that  $\inf_{t\in \R_+}\oph{2}{\bar u}{h}(t)=\kappa>0$. Then, by \eqref{bcl-eq-esti-lpt} and  \eqref{bcl-eq-esti-h}, there exist positive constants $\alpha,\beta,\eta$ such that for all $t>0$
\begin{align*}
&0<\kappa\le \nlp{h}{2}{\O}(t)\le \alpha,\\
& 0<\beta\le \oph{1}{\bar u}{u}(t)\le \eta.
\end{align*} 
 As a consequence, there exists $c_0\in\R$ such that
\begin{equation} \label{bcl-eq-F-0}
F(h(t))\to c_0 \quad \text{ and }\quad \frac{d}{dt} F(h(t))\to 0 \quad\text{ as }\quad t \to +\infty.
\end{equation}

Take now a sequence $(t_n)_{n\in\N}$ such that $t_n\to +\infty$, and consider the sequence of $L^2$ functions $(h_n)_{n\in \N}:=(h(t_n))_{n\in \N}$.  Then $\nlp{h_n}{2}{\O}$ is then bounded from above and below and therefore,   by \eqref{bcl-eq-F-0} and \eqref{bcl-cla-dF}, we get   
\begin{equation}\label{bcl-eq-dhnto0}
\lim_{n\to +\infty}\D_{2,\bar u}(h_n)= 0.
\end{equation} 
On the other hand, since $(h_n)_{n\in \N}$ is bounded in $L^2$,  there exists $\bar h\in L^2$ such that, up to extraction of a  subsequence,  $\ds{h_n\rightharpoonup \bar h}$ in $L^2$. Let us evaluate $\D_{2,\bar u}(\bar h)$. Since $\bar u$ and $m$ are positive and bounded from above and below,  the function  
$\ds{g(x):=\int_{\O}m(x,y)\frac{\bar u(y)}{\bar u(x)}\,dy}$, is well defined and $g\ge C>0$ for some positive constant $C$.  Moreover, we have  
\begin{equation}\label{bcl-eq-dhn}
 \D_{2,\bar u}(h_n)=2\left(\int_{\O}h_n^2g -\iint_{\O\times \O} m(x,y)h_n(x)h_n(y)\,dydx\right).
 \end{equation}
Since  $m\in L^{2}(\O\times \O)$ and $g\ge 0$, by Fatou's Lemma and the $L^2$ weak convergence of $h_n$, we get
\begin{align*}
&\int_{\O}\bar h^2g\le \liminf_{n\to \infty} \int_{\O}\bar h_n^2g,\\
&\lim_{n\to \infty}\iint_{\O\times \O} m(x,y)h_n(x)h_n(y)\,dydx= \iint_{\O\times \O} m(x,y)\bar h(x)\bar h(y)\,dydx.
\end{align*}

Therefore,  
$$0\le\D_{2,\bar u}(\bar h)=\iint_{\O\times \O}m(x,y)\bar u(x)\bar u(y) \left( \frac{\bar h(x)} {\bar u(x)}-\frac{\bar h(y)} {\bar u(y)}\right)^2\, dydx\le \liminf_{n\to \infty} \D_{2,\bar u}(h_n)=0,$$
which enforces that $\bar h =\nu \bar u$ for some constant $\nu \in \R$. Recall now that for all $n$, $h_n \in \bar u^{\perp}$, so
$$\nu\nlp{\bar u}{2}{\O}^2=\int_{\O}\bar h\bar u =\lim_{n\to +\infty} \int_{\O}h_n\bar u=0, $$ implying that $\nu =0$.
Now since, $\bar h=0$ and $h_n\rightharpoonup \bar h$,   from \eqref{bcl-eq-dhnto0} and \eqref{bcl-eq-dhn}, we get
$$\lim_{n\to +\infty}\int_{\O}g(x)h_n^2=0,$$ which leads to the following contradiction 

$$ \kappa\le \lim_{n\to +\infty}\int_{\O}h_n^2\le C^{-1} \lim_{n\to \infty} \int_{\O}h_n^2g=0.$$

\noindent Hence, $\inf_{t\in \R_+}\oph{2}{\bar u}{h}(t)=0$, and since $\oph{2}{\bar u}{h}(t)>0$, this implies that 
$\ds{\liminf_{t\to \infty}\oph{2}{\bar u}{h}(t)=0}.$

\fdem

\section*{Acknowledgements}

The research leading to these results has received funding from  the french ANR program under the  “ANR JCJC” project MODEVOL: ANR-13-JS01-0009 held by Gael Raoul  and the ANR "DEFI" project NONLOCAL: ANR-14-CE25-0013 held by Francois Hamel.  J. Coville also wants to thank Professor Raoul for interesting discussions on this topic.
 
%%%%%%%%%%%%%%%%%%%%%%%%%%%%%
%%%%%%%%%%%%%%%%%%%%%%%%%%%%%
%%%%%%%%%%%%%%%%%%%%%%%%%%%%%

\appendix
% 	\section{}\label{bcl-app1}
\section{Numerical Aspects}

To investigate numerically the behaviour of the solution of \eqref{bcl-eq-num}, we are led to understand how to solve numerically evolution problems of the form :
\begin{align}
&\partial_{t}v(t,x)=\int_{\O} K(x,y)v(t,y) \, dy +a(x)v(t,x)  \quad\text{ in }\quad (0,\infty)\times \O\label{fem-nonloc.eq.1}\\
%& \frac{\partial u}{\partial n}=0  \quad\text{ in }\quad \partial\O\times \R_+\\
&v(0,x)=v_0(x) \quad\text{ in }\quad \O\label{fem-nonloc.eq.1-ci}
\end{align} 
%where $c(x)$ is defined by 
%$$c(x):=\int_{\O} K(y,x)dy$$

To solve numerically \eqref{fem-nonloc.eq.1}--\eqref{fem-nonloc.eq.1-ci}, our approach is to rewrite the above problem in a variational form and  use a finite element method.
 Multiplying \eqref{fem-nonloc.eq.1} by  $w\in L^2(\O)$ and integrating over $\O$, we get  
\begin{equation}
\partial_{ t}\int_{\O}v(t,x)w -\int_{\O}\int_{\O}K(x,y) v(t,y)w(x) \, dydx -\int_{\O}a(x)v(t,x)w=0 \\
\end{equation}

To approximate the time derivation, we use standard Euler approximation scheme. For the space discretisation of  $v(t,x)$, we use the standard Lagrange finite elements.
% for the term 
%$$ \int_{\O}a(x)vw.$$
Let us look at  the non-local term. If we set,
$$v(t^n,x)\approx \sum_{i=1}^N v^{(n)}_jw_j$$
where $w_j$ is the $j$ element of a finite element basis. Then, for $w=w_i$, we have
$$\int_{\O}\int_{\O} K(x,y)v(t,y)w(x) \, dydx \approx \int_{\O}\int_{\O} K(x,y)\sum_{j=1}^N v^{(n)}_jw_j(y)w_i(x) \, dydx,$$
 which can be rewritten as follows
 \begin{equation} \int_{\O}\int_{\O} K(x,y)v(t,y)w(x) \, dydx \approx \sum_{j=1}^N v^{(n)}_j\left(\int_{\O}\int_{\O} K(x,y) w_j(y)w_i(x)\, dydx\right)\label{fem-nonloceq.approx1}
 \end{equation}

If we interpolate the map $x,y\mapsto K(x,y)$ on the fem basis $\{w_i(x)w_j(y)\}$ of $L^2(\O)\times L^2(\O) $, we have
\begin{equation}
K(x,y) \approx \sum_{k=1}^N\sum_{l=1}^N \K(k,l)w_k(x)w_l(y) \label{fem-nonloceq.approxk}
\end{equation}
Plugging the interpolation \eqref{fem-nonloceq.approxk} in the relation \eqref{fem-nonloceq.approx1}, we get 
 $$
  \int_{\O}\int_{\O} K(x,y)w_j(y)w_i(x) \, dydx \approx \sum_{k=1}^N\sum_{l=1}^N \left(\int_{\O}\int_{\O} \K(l,k) w_j(y)w_i(x)w_l(x)w_k(y)\, dydx\right) $$
  which rewrites
  \begin{equation}
  \int_{\O}\int_{\O} K(x,y)w_j(y)w_i(x) \, dydx  \approx \sum_{k=1}^N\sum_{l=1}^N \K(l,k)      \int_{\O} w_j(y)w_k(y)\, dy   \int_{\O} w_l(x)w_i(x)\, dx.\label{fem-nonloceq.approx3}
 \end{equation}
Now set $M$ and $\K$ to be the following  square matrices 
$$M_{ij}:=\int_{\O} w_i(y)w_j(y)\, dy\quad \text{ and }\quad  \K_{ij}:=K(x_i,y_j).$$
 Then \eqref{fem-nonloceq.approx3} can be expressed as follows: 
   \begin{equation}
  \int_{\O}\int_{\O} K(x,y)w_j(y)w_i(x) \, dydx  \approx \sum_{k=1}^N\sum_{l=1}^N M_{il} \K_{lk}M_{kj} \label{fem-nonloceq.approx4}
 \end{equation}
 
 The finite element matrix representing the integral term is then given by the multiplication of three matrices $M\K M$.
 Thus 
  \begin{equation} \int_{\O}\int_{\O} K(x,y)v(t,y)w(x) \, dydx \approx \sum_{j=1}^N v^{(n)}_j (M\K M)_{ij}=M\K M(v^{n})\label{fem-nonloceq.approx5}
 \end{equation}

With this finite element approximation of the integral terms, we implement a standard Euler semi-implicit scheme using FreeFem++ \cite{Hecht:2012} to compute numerically the solution of \eqref{bcl-eq-intro}--\eqref{bcl-eq-intro-ci}. To guarantee the convergence of the scheme used, the mesh used is an  adapted mesh composed of approximately  17000 triangles.

\bibliographystyle{amsplain}
\bibliography{bcl-modified.bib}
%\begin{thebibliography}
%{}\bibitem{FM}
%     {N. Fournier and S.  M{\'e}l{\'e}ard },
%       {\em A microscopic probabilistic description of a locally regulated},   {Ann. Appl. Probab.} {14},
%      {(2004)},
%   no. {4},
%        {1880--1919}.
%\end{thebibliography}
\end{document}